\documentclass[english]{article}

\usepackage{graphicx}

\usepackage{hyperref}
\usepackage{caption}

\usepackage{amsmath,amsfonts,amsthm,amssymb}
\usepackage{algorithm,algpseudocode}
\usepackage{multirow}
\usepackage{booktabs}
\usepackage{epstopdf}
\usepackage{enumerate}
\usepackage{xcolor}

\usepackage{url} 

\allowdisplaybreaks

\title{Optimization approaches to \textit{Wolbachia}-based biocontrol }
\vspace{16mm}

\author{ Jose L. Orozco-Gonzales$^{1}$, \ Antone Dos Santos$^{2}$, \ Helenice De Oliveira$^{2}$\footnote{Corresponding author: helenice.silva@unesp.br}, \\ Claudia P. Ferreira$^{2}$, \ Daiver Cardona-Salgado$^{3}$, \  Lilian S. Sepulveda-Salcedo$^{3}$, \ Olga Vasilieva$^{1}$ \\
 \vspace{10mm} \\
$^1$ \small Universidad del Valle, Cali 760031, Colombia \\
$^2$ \small S\~{a}o Paulo State University (UNESP), 18618-689 Botucatu, SP, Brazil \\
$^3$ \small Universidad Autonoma de Occidente, Cali 760032, Colombia \\
 \\
}

\date{December 22, 2023}

\begin{document}

\maketitle

\begin{abstract}
This paper proposes two realistic and biologically viable methods for \textit{Wolbachia}-based biocontrol of \textit{Aedes aegypti} mosquitoes, assuming imperfect maternal transmission of the \textit{Wolbachia} bacterium, incomplete cytoplasmic incompatibility, and direct loss of \textit{Wolbachia} infection caused by thermal stress. Both methods are based on optimization approaches and allow for the stable persistence of \textit{Wolbachia}-infected mosquitoes in the wild \textit{Ae. aegypti} populations in a minimum time and using the smallest quantity of \textit{Wolbachia}-carrying insects to release. The first method stems from the continuous-time optimal release strategy, which is further transformed into a sequence of suboptimal impulses mimicking instantaneous releases of \textit{Wolbachia}-carrying insects. The second method constitutes a novel alternative to all existing techniques aimed at the design of release strategies. It is developed using metaheuristics ($\epsilon$-constraint method combined with the genetic algorithm) and directly produces a discrete sequence of decisions, where each element represents the quantity of \textit{Wolbachia}-carrying mosquitoes to be released instantaneously and only once per a specified time unit. It turns out that a direct discrete-time optimization (second method) renders better quantifiable results compared to transforming a continuous-time optimal release function into a sequence of suboptimal impulses (first method). As an illustration, we provide examples of daily, weekly, and fortnightly release strategies designed by both methods for two \textit{Wolbachia} strains, \textit{w}Mel and \textit{w}MelPop.
\end{abstract}

\section{Introduction}
\label{}

\textit{Wolbachia} is an intracellular endosymbiotic bacterium that is not naturally present in wild \textit{Aedes aegypti} mosquitoes but can be transmitted to mosquito eggs through microinjections. When a female mosquito harbors \textit{Wolbachia} in her cells, her ability to transmit dengue and other arboviral infections when biting humans gets drastically reduced  \cite{Moreira2009,Walker2011,Aliota2016,Pereira2018}. Additionally, the bacterium provides a reproductive advantage to infected mosquitoes over non-infected insects through maternal transmission (MT) and the cytoplasmic incompatibility (CI) phenotype. These two characteristics can be perfect or imperfect depending on the ambient conditions (such as temperature, humidity, daylight length, etc.). The latter determines whether a complete population replacement or a stable coexistence of wild and \textit{Wolbachia}-infected populations can be ultimately achieved.

Both theoretical and practical approaches to the  biological control of human arbovirus infections  based on \textit{Wolbachia} have been extensively studied by numerous scholars (see the state-of-the-art reviews \cite{Caragata2021,Dorigatti2018,Ogunlade2021} and references therein).  The most common approach involves seeking strategies to release the mosquitoes carrying \textit{Wolbachia}, with the ultimate goal of achieving a long-term persistence of  \textit{Wolbachia}-carrying insects in target geographical areas. In this context, several questions naturally arise about how to conduct the releases of \textit{Wolbachia}-infected mosquitoes in practice and what aspects to consider. In light of these questions, the release programs should be biologically viable and aim to optimize economic resources, considering the cost of mosquito breeding and logistical expenses associated with releases.

The practical and theoretical challenges related to \textit{Wolbachia}-based biocontrol of mosquito populations have prompted interdisciplinary research efforts, resulting in the creation or modification of mathematical models based on differential equations. Given the initial states of infected and non-infected subpopulations of mosquitoes, these models seek to predict the population's evolution over a specific time horizon.

However, most models describing \textit{Wolbachia} propagation assumed perfect MT and CI, so the release strategies of \textit{Wolbachia}-infected insects developed based on such models may only be applicable under ideal climatic conditions. Nonetheless, several models did address not only the imperfect MT with incomplete CI \cite{Benedito2020,Farkas2017,Ferreira2020,Li2019}, but also the loss of \textit{Wolbachia} infection caused by thermal stress \cite{Adekunle2019,Ogunlade2020,Lopes2023,Orozco2023}. Therefore, a more realistic predictive model featuring imperfect MT, incomplete CI, and the loss of \textit{Wolbachia} infection due to thermal stress should be used to develop release strategies apt for a broader range of climatic conditions. In this study, we employ a comprehensive modeling framework proposed in \cite{Orozco2023} that is easily adjustable to different strains of \textit{Wolbachia}.

Different biological control strategies aimed at  partial or total population replacement of the wild mosquito population with \textit{Wolbachia}-infected mosquitoes can be tested based on a predictive model, and mathematical optimization tools play an essential role in designing such strategies. From the economic and common-sense standpoints, the best strategies should be biologically viable and require a minimum intervention time with a smaller quantity of \textit{Wolbachia}-carrying mosquitoes to release. Decisions must also be made regarding whether to conduct a single abundant release or multiple smaller-size periodic releases and determine the frequency of such releases (daily, weekly, fortnightly, monthly, etc.).

It is important to note that several publications have addressed the design of optimal \textit{Wolbachia}-based release strategies from the standpoint of dynamic optimization in continuous time and proposed continuous-time release strategies \cite{Almeida2019b,Almeida2019a,Almeida2023,Campo2017a,Campo2017b,Campo2018,Cardona2020,Cardona2021,Rafikov2019}. Even though a continuous-time function defining the optimal release strategy may yield a general structure of the release program and bring forward some other helpful insights for practitioners, there is no way to implement such a strategy in practice. In this context, impulsive and discrete control inputs seem to be more natural mathematical tools for mimicking the releases of \textit{Wolbachia}-carrying mosquitoes. Indeed, recent works \cite{Bliman2023,Vicencio2023} introduced the modeling of the releases with periodic impulses but did not address the optimization issue;  only constant release sizes of \textit{Wolbachia}-infected mosquitoes were considered.

The present study intends to go further in this strand of research by adding the optimization outlook to the design of impulsive and discrete release strategies that seek to achieve the establishment and durable persistence of \textit{Wolbachia} in the wild populations of mosquitoes in minimum time while also using the least possible number of \textit{Wolbachia}-infected insects. In particular, we propose two practical approaches for optimizing the releases of \textit{Wolbachia}-carrying mosquitoes, one adapted from the continuous-time dynamic optimization and another developed using metaheuristics directly in discrete time.

The first approach renders impulsive releases with periodic frequencies (e.g., daily, weekly, etc.) and variable release sizes, usually monotone decreasing. This method stems from dynamic optimization in continuous time and further construction of sequences of suboptimal periodic impulses. Here, we used insights from the methodology developed in \cite{Bliman2020-rus} to optimize the Sterile Insect Technique (SIT). The second approach is a novel alternative to all existing techniques aimed at the design of release strategies. It directly produces a discrete sequence of decisions, where each element represents the quantity of \textit{Wolbachia}-carrying mosquitoes to be released simultaneously and only once per a specific time unit (e.g., day, week, etc.), and the resulting sequence of decisions constitutes the release program. The second approach is formalized by the $\epsilon$-constraint method in combination with the genetic algorithm \cite{Lagaros2023} and thus renders a near-optimal solution in the form of a minimum number of effective releases and the smallest overall quantity of \textit{Wolbachia}-carrying mosquitoes to release during the intervention. As the minimal time unit is assumed to be one day, this method also allows the optimization of the day of release within more extended periods, e.g., within a week.

Both approaches are suitable for practical implementation despite the dissimilarity in the release programs they have ultimately produced. Notably, we found no references addressing the design of near-optimal discrete release strategies by genetic algorithm, and this part of our work represents a genuine contribution.

The paper is organized as follows. Section \ref{sec_modelling} presents a brief outline of the population dynamics model describing the temporary evolution of two mosquito populations sharing the same geographical space and also provides the model's core properties. This model is further used as a dynamic constraint in the optimization frameworks. Section \ref{sec-opt} covers, from the theoretical standpoint, diverse optimization approaches to design a release program facilitating the invasion of \textit{Wolbachia} with minimum time and lower intervention costs. In particular, Subsection \ref{pontryagin_results} focuses on the traditional approach based on the Pontryagin maximum principle to construct a continuous optimal strategy depicting a general structure of the optimal release program but yielding little practical usage. Nonetheless, this continuous-time optimal release strategy becomes the basis for building a sequence of suboptimal impulses that mimic the field releases well and thus exhibit practical utility (Subsection \ref{sec:subopt}). In turn, Subsection \ref{GA} gives out the main contribution of the present work, which is the construction of a discrete sequence of near-optimal releases by the genetic algorithm. Notably, this discrete sequence turns out to be very apt for practical implementation. Section \ref{results} presents the outcomes of optimization techniques obtained for two strains of \textit{Wolbachia} (\textit{w}Mel and \textit{w}MelPop) along with the core indicators of all designed release strategies. Finally, Section \ref{discussion} outlines the main results of the present study and highlights the advantages of the novel approach based on the genetic algorithm.

\section{Dynamics of the \textit{Wolbachia} invasion}
\label{sec_modelling}

In this section, we briefly present an ODE model that was proposed for describing the temporal evolution of two populations of \textit{Wolbachia}-infected and wild \textit{Aedes aegypti} mosquitoes sharing the same geographical space \cite{Orozco2023}. This model will be further laid into the basis of the dynamic optimization approaches targeting to design the optimal strategies for establishing \textit{Wolbachia} in wild mosquito populations.

\subsection{Outline of the model}

We recall first that \textit{Wolbachia} may invade a wild population of \textit{Aedes aegypti} mosquitoes by the following two reproductive mechanisms \cite{Jiggins2017}:
\begin{description}
  \item[Maternal Transmission (MT).] This mechanism allows passing the bacterium from a \textit{Wolbachia}-carrying female mosquito to all her offspring with the probability $\nu \in [0,1]$. The latter ensures that the bacterium will be consistently present in the next generation of mosquitoes, allowing it to persist.
  \item[Cytoplasmic Incompatibility (CI).] This  mechanism enables wild females to produce inviable offspring with the probability $\eta \in [0,1]$ after mating with a \textit{Wolbachia}-carrying male, thereby reducing the number of \textit{Wolbachia}-free insects in the next generation.
\end{description}

Under ideal climatic conditions, with daily variations of temperature between $20^{\circ}$C and $30^{\circ}$C, as well as in laboratory studies under controlled temperature and humidity, and depending on the bacterial strain, \textit{Wolbachia}-infected mosquitoes usually exhibit a high-level or even a perfect MT together with complete CI that is, $\nu=1, \eta=1$ \cite{Dorigatti2018}. However, the field conditions with ambient temperatures exceeding $32-35^{\circ}$C may partially disrupt both reproductive mechanisms and induce the so-called imperfect MT and incomplete CI \cite{Ross2019,Ross2017,Ulrich2016}. Consequently, the probabilities $\nu$ and $\eta$ may drop below $1$. From the biological standpoint, it would imply that a share of the offspring produced by a \textit{Wolbachia}-carrying female mosquito will be bacterium-free (with a probability $(1-\nu)$) and that a part of the offspring resulting from mating between a \textit{Wolbachia}-carrying male and a bacterium-free female will be viable and \textit{Wolbachia}-free (with a probability $(1-\eta)$).

Moreover, the intracellular density of \textit{Wolbachia} in mosquitoes may be dramatically reduced by their prolonged exposure to daily temperatures exceeding $37^{\circ}$C and thus lead to the complete loss of infection in adult insects \cite{Ogunlade2021,Caragata2023}. Following the idea of Adekunle \textit{et al.} \cite{Adekunle2019}, this detrimental effect can be modeled by the parameter $\omega \geq 0$ that denotes the per capita loss of \textit{Wolbachia} infection in adult mosquitoes due to thermal stress.

The model proposed in \cite{Orozco2023} captures the deviations in MT and CI, as well as the direct loss of \textit{Wolbachia} infection in adult insects described by the set of parameters $\big( \nu, \eta, \omega \big)$ introduced above.
It considers two populations $x(t)$ and $y(t)$ that denote, respectively, the sizes of wild and \textit{Wolbachia}-infected mosquitoes
\begin{subequations}
\label{system}
\begin{align}
\label{sys-x}
\frac{d x}{d t} &= \left( \rho_n x \frac{x + (1-\eta) y}{x + y} + (1-\nu) \rho_w  y \right) e^{-\sigma(x+y)} + \omega y - \delta_n x, \\
\label{sys-y}
\frac{d y}{d t} &= \nu \rho_w y e^{-\sigma(x+y)} - \omega y - \delta_w y + u(t),
\end{align}
\end{subequations}
where $\rho_n, \rho_w$ and $\delta_n, \delta_w$ denote the average fecundities and death rates of the wild and \textit{Wolbachia}-carrying populations, respectively, $\sigma>0$ is related to the intraspecific competition,
and $u=u(t)$ is an exogenous variable mimicking the external releases of \textit{Wolbachia}-carrying mosquitoes.

Thus, Eq. \eqref{sys-x} states that three inflows induce the growth of the wild population: (1) the total offspring of wild females, including a fraction $(1-\eta)$ of the CI survivors; (2) the fraction $(1-\nu)$ of the offspring of \textit{Wolbachia}-carrying females for which the MT fails; (3) adult insects that lose \textit{Wolbachia} infection due to thermal stress. The outflow of the wild insects is due to their natural mortality. On the other hand, Eq. \eqref{sys-y} shows that the recruitment of the \textit{Wolbachia}-carrying population is propelled by the fraction $\nu$ of the offspring of \textit{Wolbachia}-infected females for which the MT works and due to the external releases of \textit{Wolbachia}-carriers, $u(t)$, while their removal is caused not only by the natural  mortality but also by the loss of \textit{Wolbachia} infection in adult insects. Here, it is essential to highlight that all \textit{Wolbachia}-related parameters $\rho_w$, $\delta_w$, $\nu$, $\eta,$ and $\omega$ are strain-dependent and satisfy the following conditions:
\begin{equation}
\label{WB-cond}
 \rho_n > \rho_w, \qquad \delta_n < \delta_w.
 \end{equation}
The above relationships assert that wild mosquitoes are more fertile and live longer than \textit{Wolbachia}-carrying mosquitoes. This assumption is supported by scientific evidence (e.g., the state-of-art review \cite{Dorigatti2018} and references therein).

\subsection{Core properties of the model}

The ODE model \eqref{system} with $u(t)=0$ was thoroughly examined in \cite{Orozco2023}, and here we only provide a synopsis of its core properties. First, the following three basic offspring numbers are defined:
\begin{itemize}
    \item $Q_x:=\dfrac{\rho_n}{\delta_n}$, the mean number of wild individuals produced by one wild individual during its lifetime;
    \item $Q_y:=\dfrac{\nu \rho_w}{\omega+\delta_w}$, the mean number of \textit{Wolbachia}-carrying individuals produced by one \textit{Wolbachia}-infected individual during its lifetime;
    \item $Q_{y, x}:=\dfrac{(1-\nu) \rho_w+\omega Q_y}{\delta_n}$, the mean number of wild individuals produced by one \textit{Wolbachia}-carrying individual during its lifetime.
\end{itemize}
It is worth recalling that a population is viable only if its basic offspring number is above 1. Otherwise, a population is regarded as unviable, ultimately becoming extinct. Therefore, we should assume the viability of each mosquito population, $Q_x >1$ and $Q_y>1$, which, together with the conditions \eqref{WB-cond}, results in the following condition:
\begin{equation}
\label{cond-viab}
Q_x > Q_y > 1.
\end{equation}
In \cite{Orozco2023}, it was proved that, under the condition \eqref{cond-viab}, the dynamical system \eqref{system} with $u(t)=0$ is dissipative, and its absorbing set is
\[ \Omega = \left\{ (x,y) \in \mathbb{R}_+^2 : \ 0 \leq x+y \leq \dfrac{1}{\sigma} \ln Q_x  \right\}. \]
It is easy to see that all the solutions of the dynamical system \eqref{system} with $u(t)>0$ engendered by the nonnegative initial condition $\big( x(0), y(0) \big)$ will also remain nonnegative and bounded whenever the function $u(t)$ is bounded. The latter implies that for any nonnegative initial condition $\big( x(0), y(0) \big)$, the dynamical system \eqref{system} has a unique nonnegative solution for any bounded $u(t) \geq 0$, which remains bounded for all $t \geq 0$. In other words, the dynamical system \eqref{system} is well-posed and biologically meaningful for any $u(t) \geq 0$.

Second, it was shown in \cite{Orozco2023} that, under the condition \eqref{cond-viab}, the dynamical system \eqref{system} with $u(t)=0$ may possess between two and four equilibria whose existence and stability properties can be summarized as follows:
\begin{itemize}
\item The trivial $\mathbf{E}_0=(0,0)$ is a nodal repeller and the boundary equilibrium $\mathbf{E}_x=\big( x^{\sharp}, 0 \big)$, with $x^{\sharp} = \dfrac{1}{\sigma} \ln Q_x$ a nodal attractor.
\item If additionally it holds that
\begin{equation}
\label{coex-ex}
 Q_c:= \frac{Q_{y, x}+Q_y+\eta Q_x}{Q_x} >1 \quad \text{and} \quad Q_y - Q_{y,x} - 2 \sqrt{Q_{y,x} \big( Q_x - Q_y \big)} > 0,
\end{equation}
two strictly positive equilibria  $\mathbf{E}_u=(x_u,y_u)$ and  $\mathbf{E}_s=(x_s,y_s)$ arise and their coordinates
 \begin{subequations}
 \label{coex-eq}
 \begin{equation}
 \label{coex-eq-u}
        x_u = \frac{\ln Q_y}{2\eta\sigma} \left( (Q_c-1) + \sqrt{(Q_c-1)^2 - 4 \eta \dfrac{Q_{y, x}}{Q_{x}}}   \right), \quad y_u = \dfrac{1}{\sigma} \ln Q_y - x_u;
        \end{equation}
\begin{equation}
\label{coex-eq-s}
        x_s = \frac{\ln Q_y}{2\eta\sigma} \left( (Q_c-1) - \sqrt{(Q_c-1)^2 - 4 \eta \dfrac{Q_{y, x}}{Q_{x}}}   \right), \quad y_s = \dfrac{1}{\sigma} \ln Q_y - x_s.
        \end{equation}
fulfill the following relationships
\[ x_u > x_s, \qquad y_u < y_s. \]
Here, $\mathbf{E}_u$ is a saddle point, and $\mathbf{E}_s$ is a nodal attractor.
\end{subequations}
\end{itemize}
The two coexistence equilibria $\mathbf{E}_u$ and $\mathbf{E}_s$ may collide when $(Q_c-1)^2 = 4 \eta \dfrac{Q_{y, x}}{Q_{x}}$ and thus become a single saddle point. Such a situation would indicate that the dynamical system \eqref{system} with $u(t)=0$ undergoes a pitch-fork bifurcation. Furthermore, when $\nu =1, \omega=0,$ and $\dfrac{Q_x - Q_y}{Q_x} < \eta \leq 1$, the equilibrium $\mathbf{E}_s$ of stable coexistence turns into another boundary equilibrium $\mathbf{E}_y = \left( 0, \dfrac{1}{\sigma} \ln Q_y \right)$ that corresponds to the persistence of \textit{Wolbachia}-infected population and extinction of the wild mosquitoes (population replacement).

\begin{figure}[t]%
    \centering
    \begin{tabular}{ccc}
    \includegraphics[width=7cm]{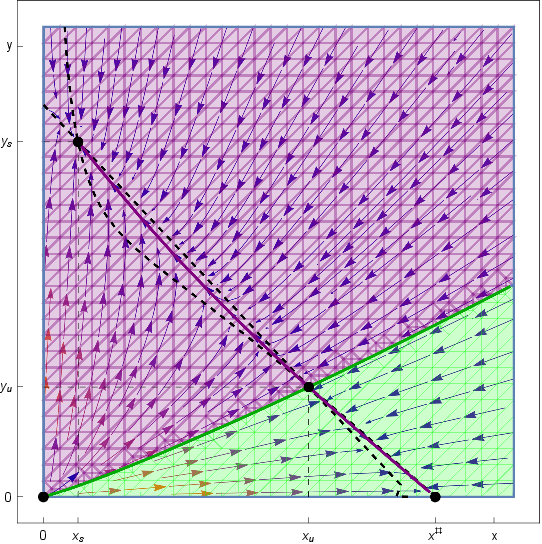} & & \includegraphics[width=7cm]{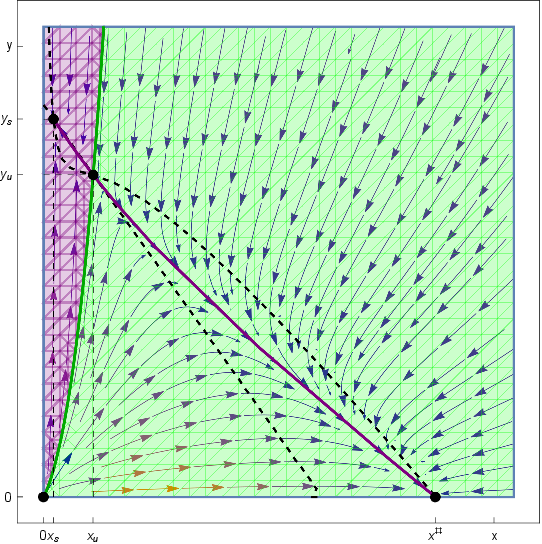} \\
    & & \\
    (a) \textit{w}Mel strain & $\qquad$ & (b) \textit{w}MelPop strain \\
    \end{tabular}
    \caption{Phase portraits to the system \eqref{system} with $u(t)=0$ under the conditions \eqref{WB-cond}, \eqref{cond-viab}, \eqref{coex-ex} and parameter values presented in Table \ref{tab:model parameters} \label{fig:phase_planes}}
\end{figure}

Thus when the conditions \eqref{WB-cond}, \eqref{cond-viab}, \eqref{coex-ex} hold, the dynamical system \eqref{system} with $u(t)=0$ exhibits bistability with two stable nodes $\mathbf{E}_x,$ $\mathbf{E}_s$ and a saddle point $\mathbf{E}_u$ in between. This situation is shown in Figure \ref{fig:phase_planes} for two \textit{Wolbachia} strains, \textit{w}Mel and \textit{w}MelPop. In both charts of Figure \ref{fig:phase_planes}, the purple solid lines depict the unstable manifold of the saddle point $\mathbf{E}_u$ that connects $\mathbf{E}_x$ and $\mathbf{E}_s$. In contrast, the green solid lines originated in $\mathbf{E}_0$ display its stable manifold that divides the positive quadrant $\mathbb{R}_+^2$ into two attraction basins. The lower basin corresponds to $\mathbf{E}_x$ (painted in light green color), and the upper basin corresponds to $\mathbf{E}_s$ (painted in light purple color). The black dashed lines on both charts of Figure \ref{fig:phase_planes} indicate the nullclines of the system \eqref{system} with $u(t)=0$.

Before the \textit{Wolbachia}-based control begins, only the wild insects are present, meaning that the initial state for the system \eqref{system} is located on the horizontal axis somewhere close to $\mathbf{E}_x$, that is, $ \big( x(0), y(0) \big)= \big( x^{\sharp}, 0 \big)$. Following the idea explored in \cite{Escobar2021}, to guarantee the evolution of the system \eqref{system} with $u(t)=0$ towards the stable coexistence equilibrium $\mathbf{E}_s$, it suffices to choose the single release size $y(0)=y_r>0$ so that $\big( x^{\sharp}, y_r \big)$ belong to the attraction basin $\mathbf{E}_s$, i.e. lay inside the purple-colored region (see Figure \ref{fig:phase_planes}). This approach may work well for the \textit{w}Mel strain of \textit{Wolbachia} (cf. Figure \ref{fig:phase_planes}(a)) because $y_r$ admits relatively low values.

However, for the \textit{w}MelPop strain of \textit{Wolbachia} (cf. Figure \ref{fig:phase_planes}(b)) the release size $y_r$ must be enormous. Alternatively, a vast single release of \textit{w}MelPop-carrying insects may be replaced by a series of periodic releases with smaller sizes. This strategy may also move the system \eqref{system} with $u(t)=0$ towards the stable coexistence equilibrium $\mathbf{E}_s$ if we add the impulsive releases of \textit{Wolbachia}-carrying insects
\[ y(0)=\Lambda, \quad y \big( i\tau^{+} \big) = y \big( i\tau^{-} \big) + \Lambda, \quad i=1,2,3, \ldots, \]
where $\Lambda>0$ stands for the release size, $\tau$ is the frequency of the releases, and $y \big( i\tau^{\pm} \big)$ denote the right and left limits of the function $y(t)$ at $t=i\tau.$ More details regarding this approach can be found in \cite{Vicencio2023}, where the interplay between the release sizes $\Lambda$ and frequencies $\tau$ was also discussed.

Even though the two works \cite{Vicencio2023,Escobar2021} provide very comprehensive ideas for the \textit{Wolbachia}-based biocontrol, they disregard the optimization aspects of this type of control intervention. Namely, the authors of these works did not focus on minimizing the total number of \textit{Wolbachia}-carrying mosquitoes to be released during the intervention nor on reducing the overall intervention time. The principal goal of the present work is to manipulate the external variable $u(t)$ to design the release strategies capable of moving the system \eqref{system} from the boundary steady state $\mathbf{E}_x$ to the attraction basin of the stable coexistence equilibrium $\mathbf{E}_s$ minimizing the time to achieve it and
the number of mosquitoes to release. It is worthwhile to point out that the external variable $u(t)$ may take different forms, such as a continuous or discrete function.

\section{Optimization approaches}
\label{sec-opt}

This section proposes several approaches for choosing the external variable $u(t)$ that denotes the number of \textit{Wolbachia}-infected mosquitoes to be released at time $t \geq 0$. These approaches are feed-forward (or have the so-called ``open-loop'' nature), meaning there is no need to assess the current sizes of both mosquito populations during the intervention, as estimating insect populations can be challenging and expensive. For all approaches considered in this section, we introduce the following magnitude constraint on the control variable $u(t)$
\begin{equation}
\label{umax}
 0 \leq u(t) \leq L,
 \end{equation}
where $L > 0$ is the maximum number of \textit{Wolbachia}-carrying insects produced (on average) per day by a facility where the transinfected mosquitoes are mass-reared. Thus, the number of \textit{Wolbachia}-carriers available to release daily is $L$, but if the releases are performed every $n$ days, the size of each release is bounded from above by the quantity $n L.$

The principal goal of the control intervention consists in moving the dynamical system \eqref{system} from its initial state $\mathbf{E}_x =\big( x^{\sharp},0 \big)$ to the attraction basin of the locally stable equilibrium $\mathbf{E}_s = \big( x_s, y_s \big)$ while minimizing the total number of released \textit{Wolbachia}-carrying mosquitoes and the overall intervention time. The releases are suspended at time $t=T^{*}$ (optimal time) when the current state of the system \eqref{system} enters ``secure region'' within the attraction basin of $\mathbf{E}_s$, i.e., satisfies the conditions
\begin{equation}
\label{basin-cond}
 x \big( T^{*} \big) < x_u, \qquad  y \big( T^{*} \big) > y_u.
\end{equation}
Here, it is worthwhile to point out that, as the wild insects possess higher individual fitness ($Q_x > Q_y$) and the density $y(t)$ of \textit{Wolbachia}-carrying insects is regulated externally via $u(t)$, the first of the two conditions in \eqref{basin-cond} is more challenging to achieve than the second one. In other words, by fulfilling the first inequality in \eqref{basin-cond}, the second inequality will also hold \cite{Campo2017a,Cardona2020}.

We first present the classical approach of continuous dynamic optimization based on the Pontryagin maximum principle (Subsection \ref{pontryagin_results}), where $u(t)$ is a continuous real function. This approach is not very convenient from the standpoint of its practical implementation. Still, it may help us emanate a general structure of the optimal release program and thus envision the planning of releases. Further, in Subsection \ref{sec:subopt}, a suboptimal strategy featuring more practicality is proposed following the idea developed in \cite{Bliman2020-rus}. Here, the continuous function $u(t)$ derived from the Pontryagin Maximum Principle is converted into a finite impulse sequence that mimics the instantaneous releases. Finally, Subsection \ref{GA} presents a heuristic approach based on the genetic algorithm that yields a discrete sequence of optimal decisions.

\subsection{Continuous optimal control}
\label{pontryagin_results}

To formalize this approach, we set the goal to define a (piecewise) continuous function $u^{*}(t) \in [0,L]$ capable of moving the dynamical system \eqref{system} from $\mathbf{E}_x$ to the attraction basin of $\mathbf{E}_s$ in a minimum time $T^{*} \in (0, \infty)$ while also minimizing the total (cumulative) number of \textit{Wolbachia}-carrying mosquitoes during the whole intervention period $[0, T^{*}]$. Here, the terminal time of the control action is free, and there are several objectives to be achieved simultaneously, namely:
\begin{enumerate}
  \item
  Moving the system  \eqref{system} from $\mathbf{E}_x$ to the attraction basin of $\mathbf{E}_s$ in a finite time.
  \item
  Minimizing the time of control intervention.
  \item
  Minimizing the total number of the released \textit{Wolbachia}-carrying mosquitoes.
\end{enumerate}

The first objective is formalized by imposing the end-point condition
\[ x(T)= x_u - 1, \]
where $x_u$ is given by formula \eqref{coex-eq-u}. This condition will ensure that $\big( x(T), y(T) \big)$ satisfies the inequalities \eqref{basin-cond}, meaning that the current state of the dynamical system is inside the attraction basin of $\mathbf{E}_s$. The second objective can be formulated in an integral form using the identity $T = \int_0^T dt$. To express the third objective, we note that the total number of released \textit{Wolbachia}-carrying mosquitoes is $\int_0^T u(t) dt$. However, we assume that the \textit{marginal} cost of \textit{Wolbachia}-based intervention is proportional to the control effort. The latter will allow us to make the objective functional more tractable from a mathematical standpoint. Thus, the second and third objectives can be fused by the following objective functional:
\[ \mathcal{J}(u,T) =  P \int \limits_0^T dt  + \dfrac{1}{2} \int \limits_0^T u^{2}(t) dt, \]
where $P>0$ denotes a weight coefficient that defines the priority of the second objective over the third one. As a result, we formulate the following optimal control problem (OCP) with free final time:
\begin{subequations}
\label{ocp}
\begin{equation}
\label{oc-func}
\min_{\begin{array}{c} 0 \leq u(t) \leq L \\ 0 < T < \infty \end{array}} \mathcal{J}(u,T)  = \int \limits_0^T \left[ P  + \dfrac{1}{2}  u^{2}(t) \right] dt
\end{equation}
subject to the dynamical system with boundary conditions
\begin{equation}
\label{oc-sys-full}
\begin{array}{lcll}
\dfrac{d x}{d t} = \left( \rho_n x \dfrac{x + (1 - \eta) y}{x + y} + (1-\nu) \rho_w  y \right) e^{-\sigma(x+y)} + \omega y - \delta_n x, & & x(0)=x^{\sharp}, & x(T)= x_u - 1, \\[5mm]
\dfrac{d y}{d t} = \nu \rho_w y e^{-\sigma(x+y)} - \omega y - \delta_w y + u(t), & & y(0)=0, &
\end{array}
\end{equation}
where the set of admissible controls is given by
\begin{equation}
\label{oc-u}
   \Big\{  u(t) \quad \text{such that} \quad 0 \leq u(t) \leq L \qquad \forall \; t \in [0,T] \Big\}.
\end{equation}
\end{subequations}

When $L$ is large enough, the OCP \eqref{ocp} has an optimal solution because the sufficient conditions given in \cite{Fleming1975} are fulfilled. Namely,
\begin{itemize}
  \item
  The set of admissible controls \eqref{oc-u} is closed, bounded, and convex in $\mathbb{R}_{+}$.
  \item
  For any bounded $u(t) \in [0, L]$, the dynamical system \eqref{oc-sys-full} has bounded solutions when $L$ is large enough.
  \item
  The right-hand side of the ODE system \eqref{oc-sys-full} is linear with respect to $u$.
  \item
  The integrand of \eqref{oc-func} is bounded from below and convex.
\end{itemize}

The OCP \eqref{ocp} can be formally solved using the variant of the Pontryagin maximum principle applicable to free terminal time problems (see, e.g., \cite{Lenhart2007} or similar textbooks) as a necessary condition for optimality. The Hamiltonian associated with the OPC \eqref{ocp} is
\begin{align}
H(\lambda_1, \lambda_2, x, y, u) = &-  P  - \dfrac{1}{2}  u^{2} + \lambda_1 \left[ \left( \rho_n x \dfrac{x + (1 - \eta) y}{x + y} + (1-\nu) \rho_w  y \right) e^{-\sigma(x+y)} + \omega y - \delta_n x \right] \notag \\[2mm]
\label{ham}
&+ \lambda_2 \: \bigg[ \nu \rho_w y e^{-\sigma(x+y)} - \omega y - \delta_w y + u \bigg],
\end{align}
where $\lambda_1$ and $ \lambda_2$ are adjoint (or co-state) variables linked to the state variables $x$ and $y$. In essence, $\lambda_1(t)$ and $ \lambda_2(t)$ are continuous real functions that express, at each $t \in [0, T]$, the marginal variations in the value of objective functional $\mathcal{J}(u, T)$ induced by the changes in current values of the state variables $x(t)$ and $y(t)$. Notably, $\lambda_1(t)$ and $ \lambda_2(t)$ are solutions to the following bidimensional ODE system
\begin{equation}
\label{sys-adj}
\dfrac{d \lambda_1}{d t}  = - \dfrac{ \partial H}{\partial x}, \qquad \dfrac{d \lambda_2}{d t}  =  - \dfrac{ \partial H}{\partial y}, \qquad \lambda_2(T)=0,
\end{equation}
and they constitute a necessary element of the Pontryagin maximum principle formulated below.

Let $\big( u^{*}, T^{*} \big)$ be optimal in the sense that
\[ \mathcal{J} \big( u^{*}, T^{*} \big) \leq \mathcal{J} (u, T) \]
for all other $u$ and $T$. Let $x^{*}(t)$ and  $y^{*}(t)$ denote the corresponding state trajectories of the system \eqref{oc-sys-full} defined for $u=u^{*}(t), t \in \big[ 0, T^{*} \big]$. Then there exist two piecewise differentiable adjoint real functions $\lambda_1(t)$ and $\lambda_2(t)$ satisfying the adjoint system \eqref{sys-adj} and the time-optimality condition
\begin{equation}
\label{time-opt}
H \big( \lambda_1(T), \lambda_2(T), x^{*}(T), y^{*}(T), u^{*}(T) \big) =0.
\end{equation}
It is well-known that the Pontryagin maximum principle permits conversion of the minimization of the objective functional \eqref{oc-func} under the constraints \eqref{oc-sys-full}-\eqref{oc-u} into a pointwise maximization of the Hamiltonian \eqref{ham} with respect to the control variable $u$ at each time $t$. As the Hamiltonian \eqref{ham} is concave in $u$, it has a critical point of maximum in $u = u^{*}$ for each fixed $t$. In other words, it holds that
\begin{equation}
\label{pmax}
H \big( \lambda_1(t), \lambda_2(t), x^{*}(t), y^{*}(t), u^{*}(t) \big) \geq H \big( \lambda_1(t), \lambda_2(t), x^{*}(t), y^{*}(t), u(t) \big)
\end{equation}
for all admissible $u(t) \in [0,L]$ and almost for all $t \in \big[ 0, T^{*} \big].$  Noting that $\dfrac{\partial H}{\partial u}= \lambda_2 - u$, the inequality \eqref{pmax} can be written in a more convenient form
\begin{equation}
\label{u-car}
u^{*}(t) = \max \Big\{  0, \min \Big\{ \lambda_2 (t), L \Big\} \Big\},
\end{equation}
usually referred to as a characterization of the optimal control, which is derived from the following relationships:
\begin{equation}
\label{pmax-cond}
\begin{array}{ccc}
u^{*}(t) = 0 & \text{if} & \dfrac{\partial H}{\partial u} < 0, \\[3mm]
0 < u^{*}(t) < L & \text{if} & \dfrac{\partial H}{\partial u} = 0, \\[3mm]
u^{*}(t) = L & \text{if} & \dfrac{\partial H}{\partial u} > 0. \\
\end{array}
\end{equation}
The optimality conditions \eqref{pmax-cond} relate the benefits and cost of the control intervention in a meaningful sense. Under the optimal release strategy $u^{*}(t)$, the marginal cost of the control action (expressed by $u^{*}(t)$ itself) matches its marginal benefit (expressed by $\lambda_2(t)$) at each $t \in \big[ 0, T^{*} \big]$. It also stems from \eqref{pmax-cond} that it is optimal either to use all available resources, i.e., set $u^{*}(t)=L$ if the marginal benefit exceeds the marginal cost (meaning that $\dfrac{\partial H}{\partial u} > 0$) or to stop the intervention by setting $u^{*}(t)=0$ (when $\dfrac{\partial H}{\partial u} < 0$).

The OCP \eqref{ocp} can be reduced to the so-called ``optimality system'' comprised by the two ODE systems, \eqref{oc-sys-full} and \eqref{sys-adj} with four boundary conditions, and with $u$ replaced by its characterization \eqref{u-car}. Thus, the optimality system is a four-dimensional boundary value problem with two initial-point conditions for $x$ and $y$ specified at $t=0$ and the other two end-point conditions for $x$ and $\lambda_2$ given at $t=T$, where $T$ is defined from the time-optimality equation \eqref{time-opt}. The latter brings forward an additional complexity that makes it virtually impossible to obtain the solution of the optimality system in a closed form. Nonetheless, a numerical solution to the optimality system and also to the OCP \eqref{ocp} can be obtained numerically, and the resulting optimal control $u^{*}(t)$ is rendered as a continuous function defined for all $t \in \big[ 0, T^{*} \big]$. Further details are provided in Subsection \ref{pontr_results}.

\subsection{Suboptimal impulsive control}
\label{sec:subopt}

Suppose that $T^{*} >0$ is the optimal time and $u^{*}(t), t \in \big[ 0, T^{*} \big]$ is the optimal control function obtained by solving numerically the OCP \eqref{ocp} using the technique presented in Subsection \ref{pontryagin_results}. The shape of $u^{*}(t)$ (approximated as a continuous real function) may certainly give us a general idea of the optimal release program. However, when it comes to the practical application, this program cannot be implemented in accurate terms. From the practical standpoint, it would be feasible to perform instantaneous periodic releases of \textit{Wolbachia}-carrying insects, instead of continuous ones, until the end-point conditions \eqref{basin-cond} are met at the optimal time $t=T^{*}$.

To plasm this idea, we propose to use a fusion of two approaches, one developed in \cite{Bliman2020-rus} for optimal implementation of the Sterile Insect Technique (SIT), and another proposed in \cite{Vicencio2023} for replacing a wild population of mosquitoes by \textit{Wolbachia}-infected ones without optimizing the underlying costs. As a result, we will construct a suboptimal release program consisting of periodic impulses bearing different frequencies and sizes.

First, let us define a formal ``extension'' $\hat{u}^{*}(t)$ of the continuous optimal control function $u^{*}(t)$
\[ \hat{u}^{*}(t) := \left\{ \begin{array}{rcl}
u^{*}(t) & \text{if} & t \in \big[ 0, T^{*} \big] \\[2mm]
0, & \text{if} & t > T^{*}
\end{array} \right. \]
and then recall that the number of \textit{Wolbachia}-carrying mosquitoes to be released during the day $n$ indicated by the optimal continuous strategy $\hat{u}^{*}(t)$ is given by
\[ U_n^{*} := \int \limits_{n-1}^{n} \hat{u}^{*}(t) \: dt, \qquad n=1,2, \ldots \hat{T}^{*}, \quad \hat{T}^{*}:=\big\lceil T^{*} \big\rceil, \]
where $\big\lceil T^{*} \big\rceil$ denotes the ceiling of $T^{*}$ if the optimal time is not an integer number and thus expresses the number of daily releases. Each element of the sequence $\Big\{ U_n^{*} \Big\} $ can be approximated using the trapezoidal rule by
\[ U_n^{tr} := \dfrac{\hat{u}^{*}(n) + \hat{u}^{*}(n-1)}{2} \approx \int \limits_{n-1}^{n} \hat{u}^{*}(t) \: dt, \qquad n=1,2, \ldots \hat{T}^{*} \]
Second, we can construct the finite sequence of suboptimal impulses $\Big\{ \widehat{U}_n^{*} \Big\} $ using the following setting:
\begin{equation}
\label{subopt-imp}
\widehat{U}_n^{*} := \left\{ \begin{array}{ccl}
\big\lceil U_n^{tr} \big\rceil,                                                 & \text{if} & U_n^{*} \leq U_n^{tr}, \\[2mm]
\max \limits_{t \in \big[ (n-1), n \big]} \big\lceil \hat{u}^{*}(t) \big\rceil, & \text{if} & U_n^{*} \geq U_n^{tr}, \end{array} \right. \qquad n= 1, 2, \ldots, \hat{T}^{*}
\end{equation}
with $\big\lceil \cdot \big\rceil$ denoting the ceilings. The elements of this sequence fulfill the condition $\widehat{U}_n^{*} \geq U_n^{*}$ for all $n= 1, 2, \ldots, \hat{T}^{*},$ and it also holds that
\[ \int \limits_{0}^{T^{*}} u^{*}(t) \: dt \leq \sum \limits_{i=1}^{\hat{T}^{*}} \widehat{U}_i^{*}, \qquad T^{*} \leq \hat{T}^{*} \]
due to which this impulsive release strategy is referred to as ``suboptimal'' because it requires a longer total time and a larger overall quantity of \textit{Wolbachia}-infected insects compared to the optimal continuous release program.

The effect of the daily suboptimal impulses $\Big\{ \widehat{U}_n^{*} \Big\} $ on the population dynamics of both mosquito species can be modeled by the following impulsive ODE system:
\begin{subequations}
\label{sys-imp}
\begin{align}
\label{sys-imp-x}
\dfrac{d x}{d t} & = \left( \rho_n x \dfrac{x + (1 - \eta) y}{x + y} + (1-\nu) \rho_w  y \right) e^{-\sigma(x+y)} + \omega y - \delta_n x, &  \hspace{-30mm} x(0)=x^{\sharp}, \\[5mm]
\label{sys-imp-y}
\dfrac{d y}{d t} & = \nu \rho_w y e^{-\sigma(x+y)} - \omega y - \delta_w y, &  \hspace{-30mm} y(0)=0, \\[5mm]
\label{sys-imp-ycond}
 y \big( n^{+} \big) & = y \big( n^{-} \big) + \widehat{U}_n^{*}, \qquad n=1,2, \ldots, \hat{T}^{*}, &
\end{align}
\end{subequations}
where $y \big( n^{\pm} \big)$ denote the right and left limits of the function $y(t)$ at $t = n, n=1,2, \ldots, \hat{T}^{*}.$

A similar rationale also allows us to construct an impulse sequence that models more sparse releases with different frequencies, e.g., every $m$ days. In this case, we define the $m$-periodic release sizes
\begin{equation}
 \label{imp-m}
 \widehat{U}_{m,i}^{*} = \sum \limits_{n=(i-1)m}^{im} \widehat{U}_n^{*}, \qquad i=1, 2, \ldots,  \hat{T}_m^{*}, \quad  \hat{T}_m^{*} := \left\lfloor \dfrac{\hat{T}^{*}}{m} \right\rfloor
 \end{equation}
where $\lfloor \cdot \rfloor$ defines the floor (i.e., the integer part) of $\dfrac{\hat{T}^{*}}{m}$ because the releases are planned for the first day of the $m-$day period. The effect of releases defined by \eqref{imp-m} on the population dynamics of both mosquito species can be modeled by the ODE system \eqref{sys-imp-x}-\eqref{sys-imp-y} with the following $m$-day periodic impulses
\begin{equation}
\label{imp-mday}
y \Big( t_{im}^{+} \Big)  = y \Big( t_{im}^{-} \Big) + \widehat{U}_{i,m}^{*}, \qquad t_{im} = 1 + (i -1)m, \qquad i=1,2, \ldots, \hat{T}^{*}_m
\end{equation}
instead of the daily impulses defined by \eqref{sys-imp-ycond}, where $y \Big( t_{im}^{\pm} \Big)$ denote the right and left limits of the function $y(t)$ at $t=t_{im}$.

However, if the frequency of the releases is comparable with the lifespan of the released insects, the conditions \eqref{basin-cond} may not be reached while using the release sizes defined by \eqref{imp-m}. In such a situation, it is recommended to construct instead of \eqref{imp-m} another $m$-periodic sequence of the release sizes using the excess estimates, that is,
\begin{equation}
 \label{imp-max}
 \widetilde{U}_{m,i}^{*} = m \times \Big\lceil \max \limits_{t \in [(i-1)m,im]} \hat{u}^{*}(t) \Big\rceil, \qquad i=1, 2, \ldots,  \hat{T}_m^{*}, \quad  \hat{T}_m^{*} := \left\lfloor \dfrac{\hat{T}^{*}}{m} \right\rfloor,
\end{equation}
to ensure the fulfillment of the conditions \eqref{basin-cond}. The effect of releases defined by \eqref{imp-max} on the population dynamics of both mosquito species can be modeled by the ODE system \eqref{sys-imp-x}-\eqref{sys-imp-y} with the following $m$-day periodic impulses
\begin{equation}
\label{imp-mday-max}
y \Big( t_{im}^{+} \Big)  = y \Big( t_{im}^{-} \Big) + \widetilde{U}_{i,m}^{*}, \qquad t_{im} = 1 + (i -1)m, \qquad i=1,2, \ldots, \hat{T}^{*}_m,
\end{equation}
where $y \Big( t_{im}^{\pm} \Big)$ denote the right and left limits of the function $y(t)$ at $t=t_{im}$. Notably, the total release sizes fulfill the following relationship
\[ \int \limits_{0}^{T^{*}} u^{*}(t) \: dt \leq \sum \limits_{i=1}^{\hat{T}^{*}_m} \widehat{U}_{m,i}^{*} \leq \sum \limits_{i=1}^{\hat{T}^{*}_m} \widetilde{U}_{m,i}^{*}, \]
where $\hat{T}^{*}_m \in \mathbb{N}$ determines the number of effective releases. Subsection \ref{pontr_results} provides examples of suboptimal impulsive release strategies with different frequencies $m>1$ where the impulse sequences are constructed either by the rule \eqref{imp-m} or \eqref{imp-max} depending on the \textit{Wolbachia} strain.

\subsection{Discrete control}
\label{GA}

Here, we propose another approach aiming to determine a discrete sequence of optimal decisions $\big( u^*(0), u^*(1), \ldots,$ $u^*(T^{*}) \big)$ to be taken at each day $t=0,1, \ldots T^{*}$ during the time horizon $[0,T^*]$ for which the final time conditions \eqref{basin-cond} on the state variables are fulfilled. Here each daily release $u^*(t), t=0,1, \ldots T^{*}$ is subject to the magnitude constraint \eqref{umax}, while the total release size $\sum \limits_{t=0}^{T^*} u^*(t)$ is minimized. Further, the minimum time $T^{*}$ satisfies the conditions  $x^*(T^*) < x_u$ and $ y^*(T^*) > y_u$, which guarantee the evolution of the dynamical system \eqref{system} to the equilibrium of stable coexistence $\mathbf{E}_s=\big( x_s, y_s \big)$ under the sequence of decisions $\big( u^*(0), u^*(1), \ldots, u^*(T^{*}) \big)$.

Thus, we encounter a multi-objective optimization problem with constraints given by the dynamical system \eqref{system} and the inequalities \eqref{basin-cond}. One goal is to reduce the time required for implementing the releases of \textit{Wolbachia}-carrying mosquitoes while satisfying the imposed constraints. The other is to minimize the total number of released insects. To explore more systematically the trade-offs between these two objectives, the impact of diverse discrete decision sequences  $\big( u(0), u(1), \ldots, u(T) \big)$ on the population dynamics of both mosquito species should be assessed and random changes in the decision sequence may entail better results. Therefore, we propose to employ a metaheuristic approach based on the genetic algorithm to test numerous combinations of daily release sizes and then identify the best sequence of discrete near-optimal decisions.

To formalize this approach, we use the $\epsilon-$constraint method (see, e.g., the recent state-of-the-art review on this topic \cite{Lagaros2023} or similar works) that involves introducing a constraint, represented by the parameter $\epsilon$, on one of the objectives while optimizing the other. For example, let us limit the final time $T$  by a value $\epsilon$, that is, $T \leq \epsilon$. Then a sequence of discrete decisions  $\big( u(0), u(1), \ldots, u(T) \big)$ satisfying the constraints  \eqref{umax}, \eqref{basin-cond} is constructed. Further, the value of $\epsilon$ is reduced, and the process is repeated. Thus, we have to solve a new sub-problem for each given $\epsilon$.

Let $S_u= \big( x(0),x(1),\ldots,x(T),y(0),y(1),\ldots,y(T) \big)$ be a solution to the ODE system \eqref{system}  with a given sequence of decisions $u= \big( u(0), u(1), \ldots, u(T) \big)$, and let $S$ denote the set of all solutions $S_u$ obtained for all possible $u$. Here, $u(t) \in \{ 0,1, \ldots, L \}$ and $t \in \{0,1, \ldots T \}$ may only take integer values. Knowing values for $\epsilon$, $x_u$, $y_u$, $x(0)$, $y(0),$ and $L$, we now propose an optimization model with decision variable $u(t)$ and state variables $x(t), y(t)$:
\begin{subequations}
\label{ga-problem}
 \begin{alignat}{1}
\label{functional}
& \min \; J (u,T) =\sum \limits_{t=1}^{T} u(t)\,\\[5pt]
\label{time_epsilon}
& T \le \epsilon  \\
\label{x_threshold}
& x(T) < x_u \\
\label{y_threshold}
& y(T) > y_u \\
\label{solution}
& S_u \in S  \\
\label{u_threshold}
& 0 \le u(t) \le L\\
\label{timespan}
& t=1,\ldots,T.
\end{alignat}
\end{subequations}
Note that the optimization problem \eqref{ga-problem} has a single objective $J$ for each given value of $\epsilon,$ and we propose a metaheuristic approach based on Genetic Algorithm (GA) to determine a near-optimal solution to the model \eqref{ga-problem}.

In GA context, each solution $i$, $u_i=\big[ u_i(0),u_i(1),\ldots, u_i(T) \big]$, to the considered optimization problem is understood as an individual and a set of solutions (individuals) as a population. The components of $u_i(t), t=0,1, \ldots,T$ of the solution are called ``genes'' and are characterized by the assigned values. To avoid possible confusion with the populations of mosquitoes $x(t)$ and $y(t)$ consisting of the insect individuals, we will employ in the sequel the following terminology when referring to the genetic algorithm. A solution to the optimization problem \eqref{ga-problem} will be called ``GA-individual,'' and the set of GA-individuals will be called ``GA-population.''  Each GA-individual can be evaluated according to its fitness, which can be measured by a function associated with the objective function of the optimization problem in question. The higher the value of the fitness function, the better the  GA-individual.

The GA modifies a GA-population consisting of $N$ GA-individuals (fixed size) through the genetic operators known as Selection, Crossover and Mutation. The GA-population then evolves toward a  near-optimal solution over successive iterations (called GA-generations). The Selection is generally stochastic, and the GA-individuals  with higher fitness function values are more likely to be selected. This method is the most common and is called ``Tournament Selection'' \cite{Lagaros2023} because it mimics a tournament-style competition among individuals. The winner's genes are selected to create offspring, i.e., to undergo Crossover and Mutation.

\begin{figure}[h]%
    \centering
   \includegraphics[width=11cm]{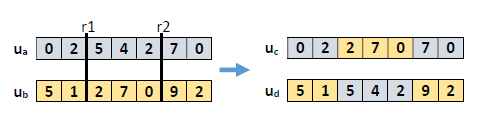}
\caption{Illustration of the Crossover operation with two crossover points $r1=2$ and $r2=5$.}
\label{cros}
\end{figure}

Crossover operation mimics the process of genetic recombination in natural evolution. The core idea is to exchange some traits of two existent GA-individuals from the pool (i.e., ``parent'' GA-individuals) and thus create one or more new GA-individuals (i.e.,  ``offspring'' GA-individuals) that are added to the GA-population. Figure \ref{cros} illustrates the Crossover operation with two designated crossover points, which is equivalent to performing two single-point crossovers. Here, the original GA-individuals $u_a$ and $u_b$ are parents, and the new GA-individuals $u_c$ and $u_d$ represent the offspring.

\begin{figure}[h]%
    \centering
   {{\includegraphics[width=5.5cm]{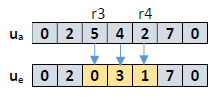} }}
\caption{Illustration of the Mutation operation on the gene's sequence from gene $r3=3$ to gene $r4=5$.}
\label{multa}
\end{figure}

Mutation operation is analogous to the biological process of genetic mutation, where small random changes occur in one or more genes. The mutation is crucial for introducing and maintaining genetic diversity in the GA-population, and it prevents premature convergence of the algorithm to a locally optimal solution. The probability of each gene undergoing a mutation (called ``mutation rate'')  is generally assumed to be around $0.05$. However, higher mutation rate values may allow for more significant exploration of the solution space while still focusing on known high-quality regions to refine and improve solutions.  Figure \ref{multa} illustrates a mutation of the GA-individual $u_a$, generating the new GA-individual $u_e$.

Thus, the GA-population will comprise more than $N$ GA-individuals at the end of each iteration due to crossover and mutation. Therefore, a new GA-population must be created by retaining $N$ best-fitted GA-individuals and discarding the others.

Another important concept related to evolutionary algorithms is elitism. In the GA-context, elitism is a strategy aimed at preserving the best GA-individuals through generations. It consists of retaining the  best-quality solutions (GA-individuals in each GA-population) and improving convergence to a globally optimal solution. A certain number $M<N$ (or a percentage) of the best-fitted GA-individuals from the current GA-generation are identified and directly copied to the next GA-generation without undergoing Crossover or Mutation operations. Further, these $M$ elite GA-individuals are updated through the best fitness of each GA-generation.

Terminating the search performed by evolutionary algorithms is a crucial task, and it is pervasive to use a combination of stopping criteria to enhance the efficacy of the termination strategy and computational efficiency while meeting the characteristics of the available computational resources. Here are the standard stop criteria for GAs that can also be combined:
\begin{itemize}
  \item
  \textit{A satisfactory solution with sufficient fitness is found. } The algorithm terminates when a solution meeting a satisfactory fitness level is found, even if it is not necessarily the optimal solution. This criterion is applicable when the goal is to find a solution that meets certain criteria rather than optimizing indefinitely.
  \item
  \textit{Stagnation or plateau is detected.} The algorithm terminates if there are no significant changes or improvements in the GA-population over a specified number of GA-generations (stagnation) or if it detects that the current GA-population has reached a ``plateau`'', where all GA-individuals are very similar and further exploration may not lead to significant improvements.
  \item
  \textit{Convergence criteria.} The algorithm terminates if almost no improvement is observed over a certain number of consecutive GA-generations. This might be measured by a ``convergence index'' that quantifies the changes in the average fitness or the best fitness in the GA-populations through GA-generations. If such a change falls below a threshold defined by the convergence index, it may indicate that further iterations are unlikely to yield significant improvements.
  \item
  \textit{Resource limits and user-defined criteria.} The algorithm terminates after a predefined number of GA-generations have been created and evaluated, a maximum number of function evaluations has been performed, or when the CPU time limit or memory usage is reached.
\end{itemize}

The Algorithm \ref{alg:cap} illustrates the GA realization proposed for solving the optimization problem \eqref{ga-problem} with different frequencies of releases $p$: $p =1$ (day), $p=7$ (week), and $p=14$ (fortnight, or two weeks). Here, the releases of \textit{Wolbachia}-carrying mosquitoes may occur every day, once per week, or once per fortnight. In the description of Algorithm \ref{alg:cap}, $N$ defines the GA-population size (number of GA-individuals in each GA-generation), $M=1$ is the size of elite population, and $G$ stands for the maximum number of GA-generations (or the maximum number of iterations to be carried out). At the beginning of each iteration, the  GA-population is denoted by the set $P$ that further grows to a larger set $P_a \supset P$ after the Crossover and Mutation operations are performed.

Each GA-individual $i$ in the population $P$ or $P_a$ is a control strategy $u_i$ to be applied in a horizon time $[1,T]$  with the following structure:
\[ \big[ u_i(1)\quad	u_i(2)\quad	u_i(3)\quad	\dots \quad	u_i(T) \big], \]
where the element $u_i(t),\;  t=1,\dots,T$ takes integer values between $0$ and $L$ and expresses the quantity of the \textit{Wolbachia}-carrying mosquitoes to be released on the day $t$. Note that $u_i(0)=0$ is assumed for all GA-individuals since the control action starts later than $t=0$. In the case of weekly or fortnightly releases ($p=7$ or $p=14$), on each time interval $\big[  pn+1, p(n+1) \big]$, $n=0,1,2\dots$, there will be only one single $t=j \in \Big\{pn+1, \ldots, p(n+1) \Big\}$ for which $u_i(j)$ may take integer values between $0$ and $pL$, so that $u_i(t)=0, \;\forall t \neq j$.

The fitness function of each GA-individual $u_i(t)$ is defined as follows:
\begin{equation}
\label{fit-fun}
 F(u_i) = \dfrac{1}{J(u_i) + pLT \: \mathbf{I}_c(u_i)}, \quad \text{where} \quad \mathbf{I}_c(u_i) := \left\{
\begin{array}{cl}
0 & \text{if }  S_{u_i} \text{ fulfills \eqref{x_threshold}, \eqref{y_threshold}} \\[2mm]
1 & \text{otherwise} \\
\end{array} \right.
\end{equation}
is a binary indicator function that penalizes the violation of inequality constraints \eqref{x_threshold} and \eqref{y_threshold}.

\begin{algorithm}[htpb]
\caption{Genetic Algorithm}\label{alg:cap}
{\small
\begin{algorithmic}
\State{\textbf{1. Enter the initial data values} }
\State{Enter with values for all the optimization model \eqref{ga-problem} parameters associated to a specific strain, and $p$, $G$, $L$, $N$, $M$, and $T$ as a multiple of $p$.}

\State{\textbf{2. Create the initial GA-population $P$ with $N$ GA-individuals}}
\State{The initial GA-population $P$ must contain $N$  GA-individuals $u_i(t), t = 1,2, \ldots, T, i =1,2, \ldots, N$ all of which satisfy the constraints \eqref{u_threshold}.}

\State{\textbf{3. Evaluate each GA-individual from $P$}}
\State{Calculate the fitness value for each individual $i$ in $P$ using the  function $F(u_i)$ defined by \eqref{fit-fun}.}

\State{\textbf{4. Store the Elite GA-individual}} \par
$u_{elite} \gets \bar{u}_i$,  such that $F(\bar{u}_i) = \max \limits_{i=1,\ldots, N }  F(u_i)$.

\State{\textbf{5. Apply genetic operators to GA-individuals of the current GA-population $P$}} \par
\State{\textbf{5.1. Selection}}
\State{Use the Tournament Selection method:}

\For{\texttt{$k=1:N$}}
      \State{1. Determine two integer random number $r$ and $s$ in $[ 1, N]$.}
      \State{2. Select the individual whose Fitness value is greater:}
                \If{$F(u_r) > F(u_s)$}
                \State $u_{k}^{selected} \gets u_r$
                \Else
                \State $u_{k}^{selected} \gets u_s$.
                \EndIf
 \EndFor

\State{\textbf{5.2. Crossover}} \par
\State{Randomly take two pairs of the selected GA-individuals $u_{k}^{selected}$.}
\State{To each pair, apply the two-point crossover technique: choose two randomized gene positions $r1$ and $r2$ (multiples of $p$) on the selected GA-individuals and perform the crossover as shown in Figure \ref{cros}.}

\State{\textbf{5.3. Mutation}} \par
\For{\texttt{all individuals}}
\State{Randomly take one value $v \in [0,1]$}
\If{\texttt{$ v<0.05$}}
\State{Perform the mutation: }
\State{Randomly take two integer values $r3 \in [1, N-1]$ and $r4 \in [r3,N]$ (multiples of $p$).}
\State{Assign a random integer value taken in interval $[0, pL]$ to the genes from $r3$ through $r4$, if $p=1$, as shown in Figure \ref{multa},} \State{or to the non-zero genes $j$ if $p=7$ or $p=14$.}

\EndIf
\EndFor

\State{\textbf{6. Evaluate all created new GA-individuals}}
\State{Calculate the fitness function value for all newly created GA-individuals according to formula \eqref{fit-fun}.}

\State{\textbf{7. Create the new GA-population}}
	 \State{$P_a \gets \Big\{P \: \bigcup \: \{ \text{GA-individuals created by crossover} \} \: \bigcup \: \{ \text{GA-individuals created by mutation} \} \Big\} $, }
	 \State{$P \gets \{ \}$}
	 \State{$P \gets \big\{ N \text{ best individuals from } P_a \big\}$,}
	 \State{$P_a \gets \{\}$.}

\State{\textbf{8. Store the Elite GA-individual}}
\State{$u_{elite} \gets \bar{u}_i$  such that $F(\bar{u}_i)= \max \limits_{i=1,...,N} F(u_i)$.}

\State{\textbf{9. Apply Stop Criterion}}
\If{\texttt{the maximum number of generations $G$ is reached}}
\State{go to \textbf{10}}
\Else
\State{go to \textbf{5}.}
\EndIf

\State{\textbf{10. End. Take the best solution $u^*$}}
\State{$u^* \gets u_{elite}$.}
\end{algorithmic}
}
\end{algorithm}

%
%


\section{Numerical results}
\label{results}

In this section, we illustrate the optimization approaches proposed in Section \ref{sec-opt} by presenting different release strategies $u(t)$ (continuous, impulsive, and discrete) capable of driving the dynamical system \eqref{system} from its initial  \textit{Wolbachia}-free state $\mathbf{E}_x$ to the attraction basin of the stable coexistence equilibrium $\mathbf{E}_s$, while minimizing the total number of released  \textit{Wolbachia}-carrying mosquitoes and the overall intervention time.

Before we proceed, it is worthwhile to recall that, in the last decades, various strains of \textit{Wolbachia} underwent laboratory tests for their virus-blocking abilities in \textit{Aedes aegypti} mosquitoes, during which these strains also showed dissimilar levels of individual fitness loss in mosquito hosts (see, e.g., state-of-the-art review by Dorigatti, \textit{et al}. \cite{Dorigatti2018}). In particular, the \textit{w}Mel strain, which exhibits a moderate virus-blocking ability and induces a weaker fitness loss in mosquitoes successfully invaded the wild \textit{Ae. aegypti} populations after field releases \cite{Ross2022,Gesto2021,Tantowijoyo2020}. However, the \textit{w}MelPop strain, bearing the highest level of virus blockage but more damaging to mosquitoes' fitness, failed to invade the wild \textit{Ae. aegypti} population \cite{Nguyen2015} despite some efforts \cite{Yeap2014,Ritchie2015}. This previous experience clearly illustrates that the release programs must be strain-dependent.

Throughout this section, we will focus on two different strains of  \textit{Wolbachia}, \textit{w}Mel and \textit{w}MelPop, which bear contrasting characteristics related to their virus-blocking abilities and individual fitness loss in mosquitoes. Table \ref{tab:model parameters} summarizes the parameters of the system \eqref{system}, their corresponding values used in numerical simulations, and the underlying references. It is worth pointing out that the initial population of wild mosquitoes $x^{\sharp}= 7030$ (either males or females) corresponds to the average mosquito density per one hectare \cite{Bliman2019}, while the coordinates of the saddle-point unstable equilibrium $\mathbf{E}_u = \big( x_u, y_u \big)$ are strain-dependent, that is,
\[ \text{ \textit{w}Mel strain } \Rightarrow \ \mathbf{E}_u = (4592,1793); \quad \text{ \textit{w}MelPop strain } \Rightarrow \ \mathbf{E}_u = (1050,3778). \]
The above values will be further used in the boundary condition assigned to the system \eqref{oc-sys-full} and in the inequalities \eqref{x_threshold}-\eqref{y_threshold}.

\begin{table}[t]
\begin{center}
\caption{Values of the model's parameters used in numerical simulations}
\begin{tabular}{llc}
\hline
\textbf{Description} & \textbf{Assumed value} & \textbf{References} \\ \hline
Fecundity rate of uninfected insects                        & $\rho_n=4.55$ days$^{-1}$                  & \cite{Styer2007} \\
Fecundity rate of \textit{w}Mel-infected insects            & $\rho_w=0.9\times\rho_n=4.095$ days$^{-1}$ & \cite{Dorigatti2018,Hoffmann2015} \\
Fecundity rate of \textit{w}MelPop-infected insects         & $\rho_w=0.5\times\rho_n=2.275$ days$^{-1}$ & \cite{Dorigatti2018,McMeniman2009} \\
Natural mortality rate of uninfected insects                & $\delta_n=1/28=0.0357$         days$^{-1}$ & \cite{Styer2007} \\
Natural mortality rate of \textit{w}Mel-infected insects    & $\delta_w=\delta_n/0.9=0.0397$ days$^{-1}$ & \cite{Dorigatti2018,Hoffmann2015} \\
Natural mortality rate of \textit{w}MelPop-infected insects & $\delta_w=\delta_n/0.5=0.0714$ days$^{-1}$ &  \cite{Dorigatti2018,McMeniman2009} \\
Feature related to mosquito competition                     & $\sigma=0.1/140 $  individual$^{-1}$       & \cite{Bliman2019} \\
MT probability in \textit{w}Mel-infected insects            & $\nu=0.95$                     & \cite{Ogunlade2021,Hoffmann2015} \\
MT probability in \textit{w}MelPop-infected insects         & $\nu=0.99$                     & \cite{Ogunlade2021,Hoffmann2015} \\
CI probability in \textit{w}Mel-infected insects            & $\eta=0.98$                    & \cite{Ogunlade2021,Hoffmann2015} \\
CI probability in \textit{w}MelPop-infected insects         & $\eta=0.95$                    & \cite{Ogunlade2021,Hoffmann2015} \\
Loss rate of \textit{w}Mel infection                        & $\omega = 0.001$ days$^{-1}$              & \cite{Ogunlade2021,Hoffmann2015} \\
Loss rate of \textit{w}MelPop infection                     & $\omega = 0.00015$ days$^{-1}$            & \cite{Ogunlade2021,Hoffmann2015} \\
\hline
\end{tabular}
\label{tab:model parameters}
\end{center}
\end{table}

As the strains affect differently the individual fitness of mosquito host (cf. values of $\rho_w$ and $\delta_w$ in Table \ref{tab:model parameters}), the average upper daily limit $L$ sufficient for the existence of the optimal control may also differ depending on the strain. Here, we assume the following values for $L$ in all numerical simulations presented in this section:
\begin{equation}
\label{L-value}
\text{ \textit{w}Mel strain } \Rightarrow \ L=750 \text{ individuals}; \quad \text{ \textit{w}MelPop strain } \Rightarrow \ L=1000 \text{ individuals.}
\end{equation}

Under these values of $L$, the optimal control problem \eqref{ocp} and the optimization problem \eqref{ga-problem} have feasible numeric solutions.
In the sequel, we present the strain-dependent optimal release strategies that stem from different approaches. Subsection \ref{pontr_results} encompasses the optimal (continuous) and suboptimal (impulsive) strategies originating from continuous dynamic optimization. In contrast, Subsection \ref{ga_results} illustrates a discrete heuristic approach where the genetic algorithm designs the optimal release strategies. For both approaches, we consider three release frequencies, that is, daily, weekly, and fortnightly (every two weeks) releases of the insects carrying either \textit{w}Mel or \textit{w}MelPop strain of \textit{Wolbachia}.

\subsection{Release strategies derived from the Maximum Principle}
\label{pontr_results}

To solve the optimality system resulting from the optimal control problem \eqref{ocp}, we have employed the advanced solver GPOPS-II\footnote{Note that GPOPS-II Manual with basic descriptions is downloadable from \url{http://www.gpops2.com/}} (Next-Generation Optimal Control Software, developed by the Optimization Research LLC, Gainesville, FL, USA). This software package implements the numerical technique based on the direct orthogonal collocation \cite{Rao2010} by scaling an input interval $[0, T]$ (with $T$ as an initial guess) to the interval $[-1, 1]$. The latter allows the users to deal with time-optimal problems and to obtain $T^{*}$ as an output only if the time-optimality condition \eqref{time-opt} is eventually met.

In the formulation of the objective functional $\mathcal{J}(u,T)$ (see formula \eqref{oc-func}), the weight coefficient $P>0$ defined the priority of the time minimization over the minimization of the total number of  \textit{Wolbachia}-carriers to be released during the intervention. As shown in other works (see, e.g., \cite{Campo2017a,Campo2017b,Campo2018}), different values of $P>0$ result in different optimal release strategies. Moreover, in the mentioned works, the authors have identified, for a fixed the upper daily limit $L$, a tradeoff between the intervention time $T^{*}$, the total number of \textit{Wolbachia}-carrying mosquitoes to be released, $\int \limits_0^{T^{*}} u^{*}(t) dt,$ and the priorities of the decision-making expressed by $P>0$ in the formulation of $\mathcal{J}(u,T)$. Namely, the higher values of $P>0$ render a shorter overall time $T^{*}$ of the intervention and a larger number of insects to be released. Alternatively, the lower values of $P>0$ entail a longer overall time $T^{*}$ of the intervention and a smaller number of insects to be released. In our simulations, we assign the higher priority to the time optimization and assume in the sequel that $P=10^6$.

\begin{figure}[t]%
    \centering
    \begin{tabular}{ccc}
        (a) & & (b) \\ & & \\
           \includegraphics[width=7cm]{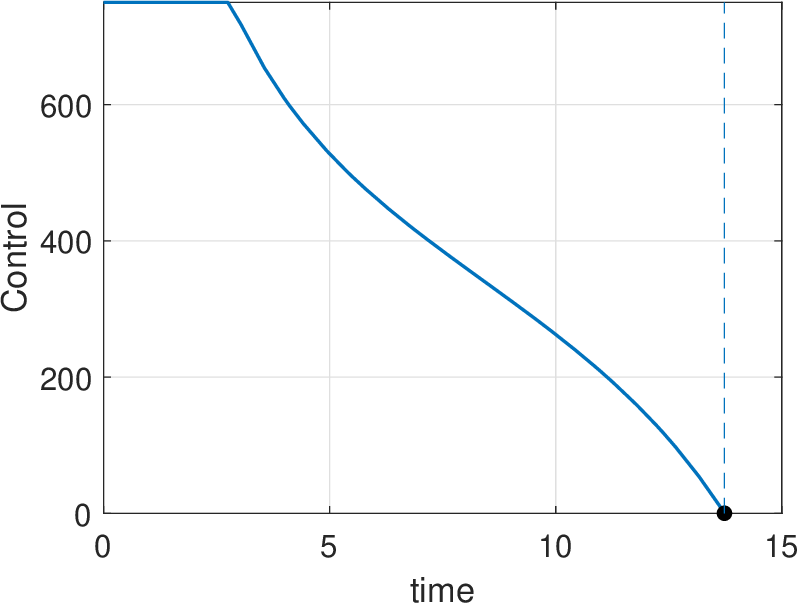} & \qquad & \includegraphics[width=7cm]{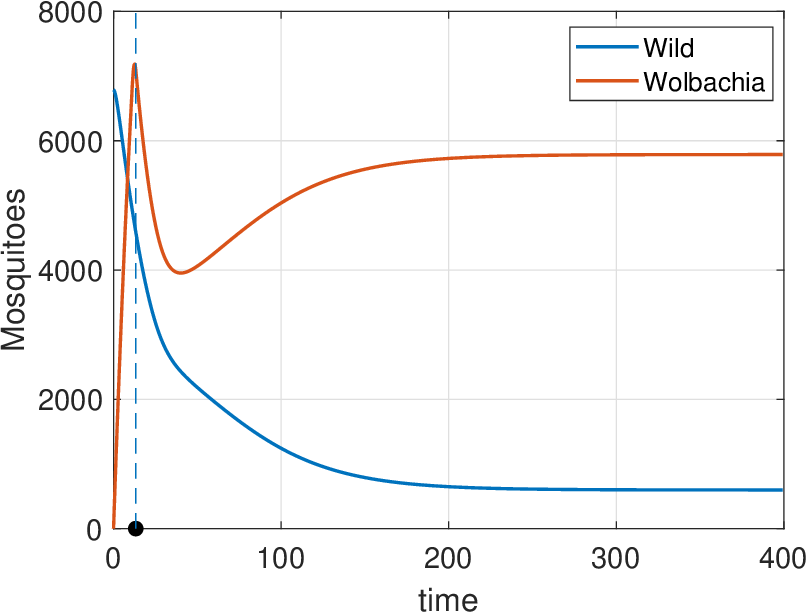} \\
    \end{tabular}
    \caption{ Numerical solutions to the OCP \eqref{ocp} for \textit{w}Mel strain: (a) optimal control $u^{*}(t), t \in \big[ 0, T^{*} \big]$ with $T^{*}=13.72$ days; (b) optimal state trajectories $x^{*}(t)$ and $y^{*}(t), \ t \in [0, 400]$  \label{figPO1}}
\end{figure}

\begin{figure}[t]%
    \centering
    \begin{tabular}{ccc}
     (a) & & (b) \\ & & \\
    \includegraphics[width=7cm]{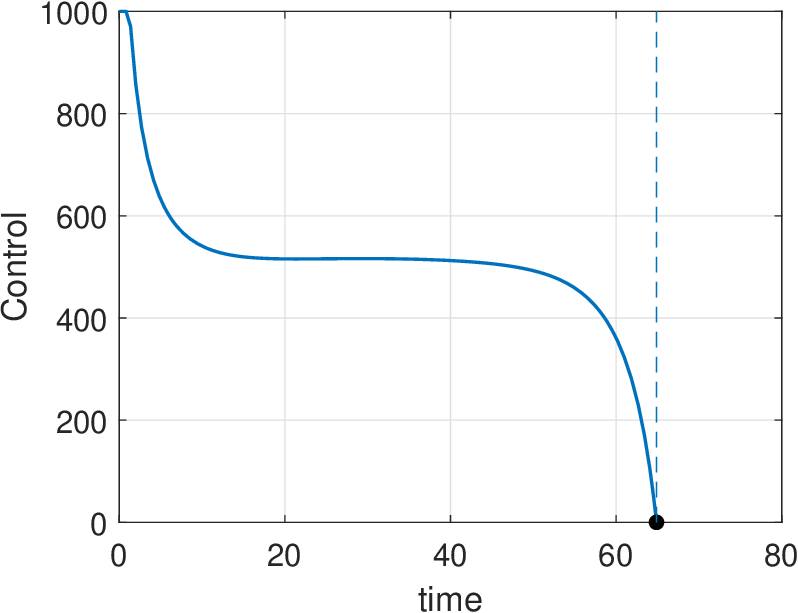} & & \includegraphics[width=7cm]{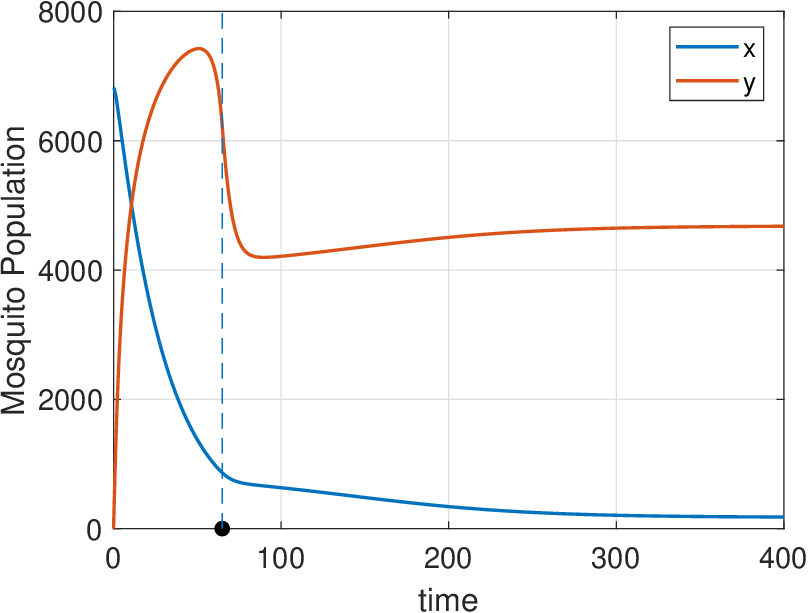} \\
    \end{tabular}
    \caption{Numerical solutions to the OCP \eqref{ocp} for \textit{w}MelPop strain: (a) optimal control $u^{*}(t), t \in \big[ 0, T^{*} \big]$ with $T^{*}=64.87$ days; (b) optimal state trajectories $x^{*}(t)$ and $ y^{*}(t), \ t \in [0, 400]$    \label{figPO2}}
\end{figure}

Numerical solutions to the OCP \eqref{ocp} for the \textit{w}Mel and \textit{w}MelPop strains of \textit{Wolbachia} are presented in Figures \ref{figPO1} and \ref{figPO2}, respectively. The left-hand charts on both figures display the optimal controls $u^{*}(t)$ as continuous functions of $t \in \big[ 0, T^{*} \big]$ that express the general structures of the optimal release programs. Comparing the left-hand charts in Figures \ref{figPO1} and \ref{figPO2}, we observe that $u^{*}(t)$ exhibit strikingly dissimilar structures for both \textit{Wolbachia} strains. From some published works, we know that the optimal release strategies based on the Pontryagin maximum principle may have different forms and structures depending on the formulation of the objective functional and other side conditions, as well as on the  \textit{Wolbachia} strain (see, e.g., \cite{Cardona2021}). Namely, such optimal strategies may exhibit a shape of a bell \cite{Almeida2019b,Campo2017a,Campo2018,Cardona2020,Cardona2021}, or a curvilinear trapezoid form \cite{Almeida2019a,Almeida2023,Campo2017a,Campo2017b,Cardona2021,Rafikov2019}, or even have a bang-bang structure \cite{Almeida2019b,Almeida2019a,Almeida2023}. Here, it is worth noting that the optimal strategy obtained for releasing the mosquitoes carrying \textit{w}Mel (cf. Figure \ref{figPO1}(a)) stays in line with the previous findings as it has a form of a curvilinear trapezoid. On the other hand, the optimal strategy obtained for releasing the mosquitoes carrying \textit{w}MelPop (cf. Figure \ref{figPO2}(a)) has an interesting form that did not arise in other studies.

Figures \ref{figPO1} and \ref{figPO2} also show that the optimal time $T^{*}$ (marked by the dashed vertical line in all four charts) is shorter for the \textit{w}Mel strain (about 14 days) compared to the \textit{w}MelPop strain (about 65 days). This outcome is rather expectable since the \textit{w}MelPop-infected mosquitoes have lower fecundity and higher mortality than \textit{w}Mel-infected ones. Moreover, the total number of \textit{Wolbachia}-carrying insects to be released during the intervention can be assessed by calculating the area below the curve $u^{*}(t), t \in \big[ 0, T^{*} \big]$ and we obtain the following estimations:
\begin{equation}
\label{int-value}
\text{ \textit{w}Mel strain } \Rightarrow \ \int \limits_0^{T^{*}} u^{*}(t) \ dt = 5961; \quad \text{ \textit{w}MelPop strain } \Rightarrow \ \int \limits_0^{T^{*}} u^{*}(t) \ dt = 33125.
\end{equation}
Thus, establishing the \textit{w}MelPop strain in a wild mosquito population is longer and more expensive compared to the \textit{w}Mel strain.

On the other hand, the right-hand charts in Figures \ref{figPO1} and \ref{figPO2} show that, after suspension of releases ($t > T^{*}$), the optimal trajectories $x^{*}(t)$ and  $y^{*}(t)$ converge to the locally stable equilibrium of coexistence $\mathbf{E}_s=\big( x_s, y_s \big)$ whose coordinates are
\[ \text{ \textit{w}Mel strain } \Rightarrow \ \mathbf{E}_s = (598,5787); \quad \text{ \textit{w}MelPop strain } \Rightarrow \ \mathbf{E}_s = (135,4693). \]
These values clearly indicate that, after \textit{Wolbachia} is locally established, the remanent size $x_s$ of the wild population (fully capable of transmitting the virus) will be more than four times larger in the case of \textit{w}Mel strain compared to the \textit{w}MelPop strain.

\begin{figure}[H]
\centering
     Daily releases \\ \vspace{4mm}
        \begin{tabular}{ccc}
    \includegraphics[width=7cm]{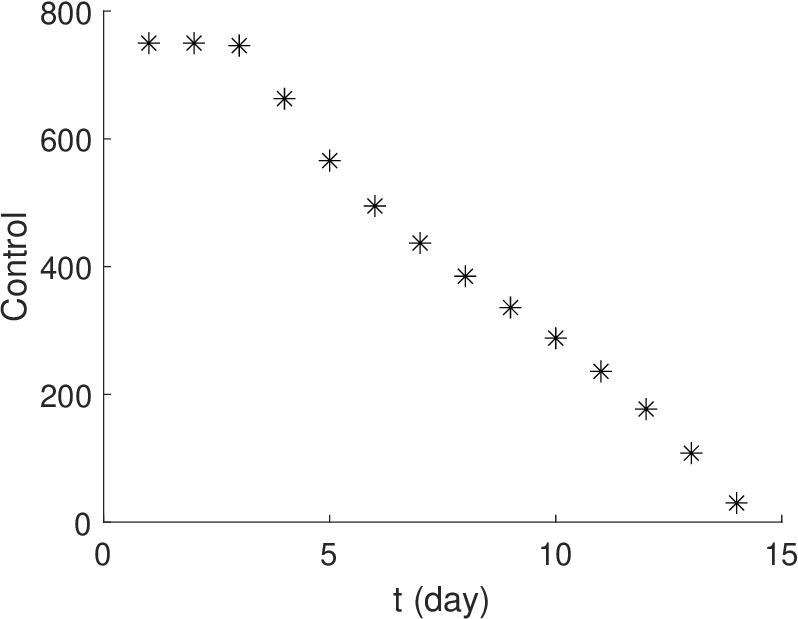} & &\includegraphics[width=7cm]{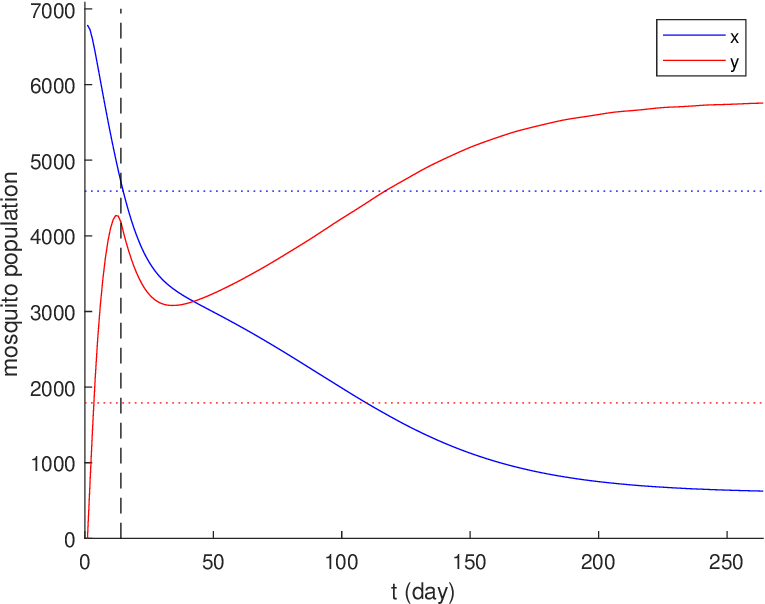} \\
    \end{tabular}
    \vspace{2mm}

    Weekly releases \\ \vspace{4mm}
         \begin{tabular}{ccc}
    \includegraphics[width=7cm]{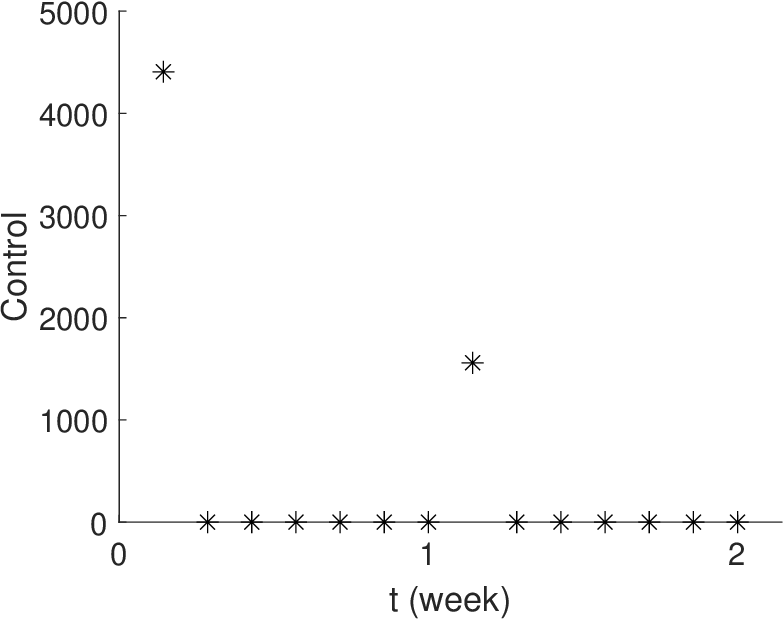} & &  \includegraphics[width=7cm]{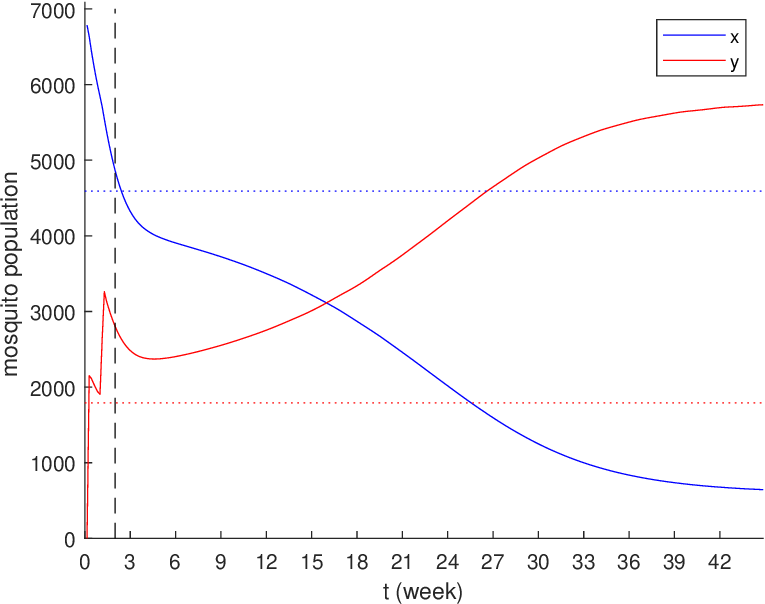} \\
    \end{tabular}
     \vspace{2mm}

    Fortnightly releases \\   \vspace{4mm}
     \begin{tabular}{ccc}
    \includegraphics[width=7cm]{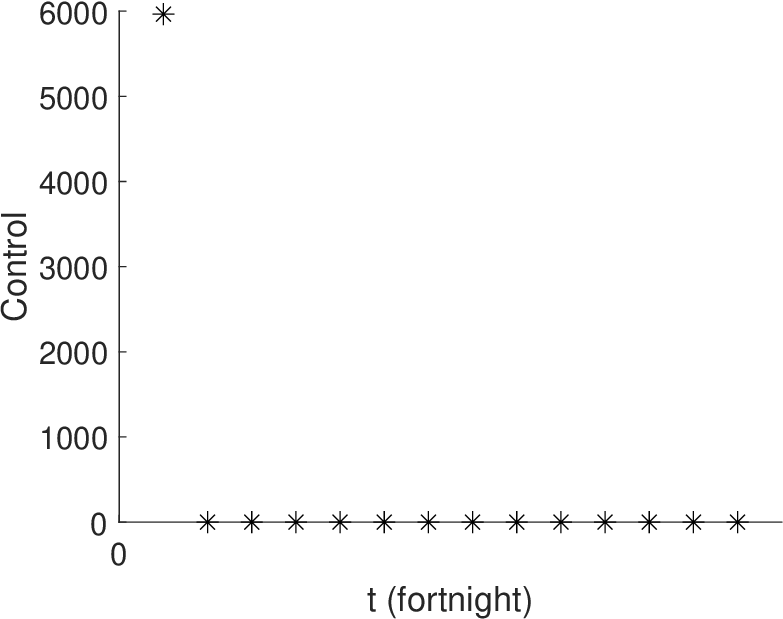} & & \includegraphics[width=7cm]{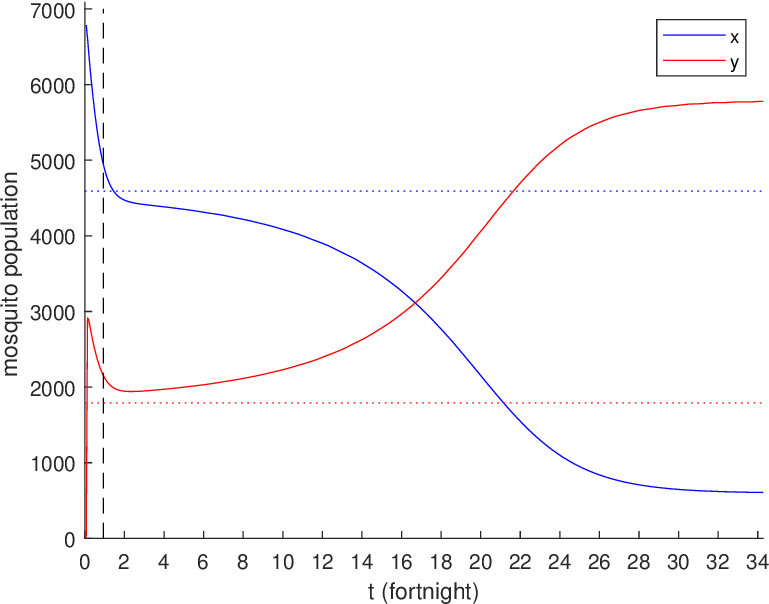} \\
    \end{tabular}
    \caption{Suboptimal impulsive control strategies for releasing \textit{w}Mel-carrying mosquitoes \textit{(left column)} and the corresponding evolution of wild mosquitos $x(t)$ and \textit{w}Mel-infected mosquitoes $y(t)$  \textit{(right column)}, where the vertical line marks the time of fulfillment of the condition \eqref{basin-cond}. \label{fig:wMel-subopt} }
\end{figure}

\begin{figure}[H]
    \centering
     Daily releases \\ \vspace{4mm}
        \begin{tabular}{ccc}
    \includegraphics[width=7cm]{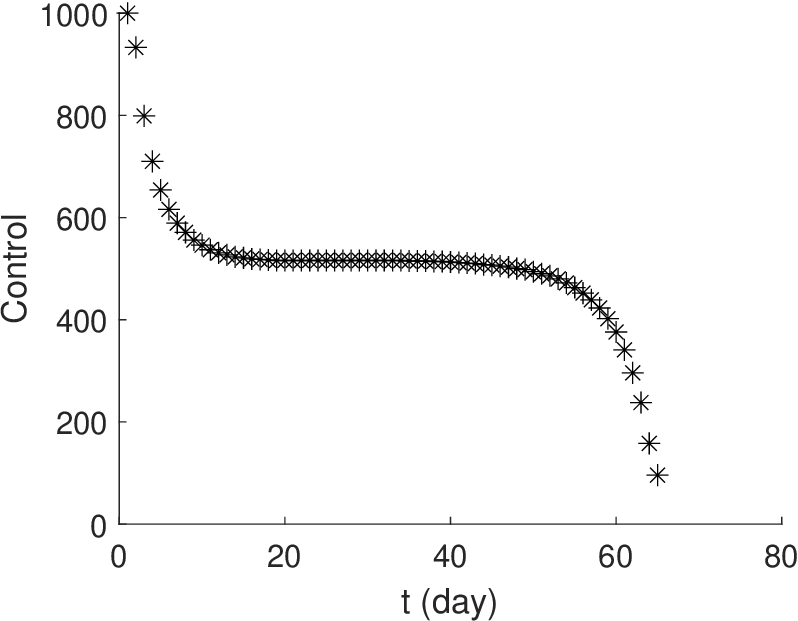} & & \includegraphics[width=7cm]{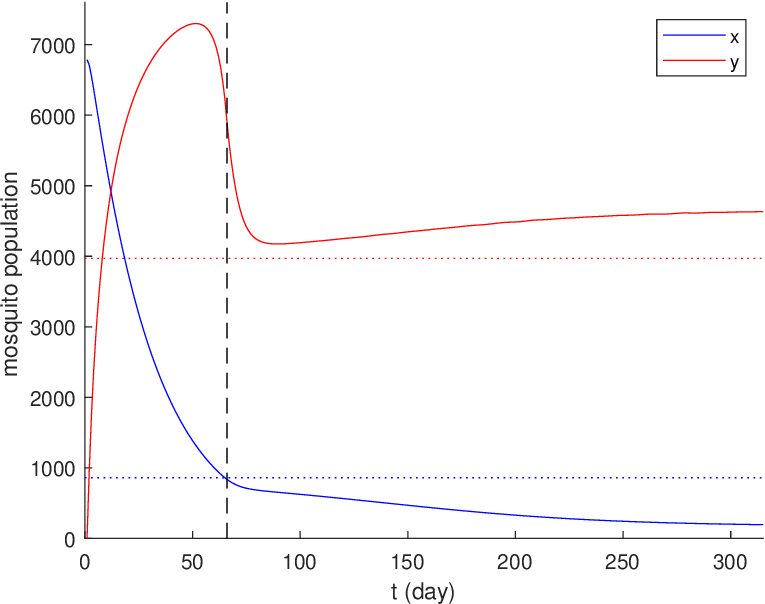} \\
    \end{tabular}
    \vspace{2mm}

    Weekly releases \\ \vspace{4mm}
         \begin{tabular}{ccc}
    \includegraphics[width=7cm]{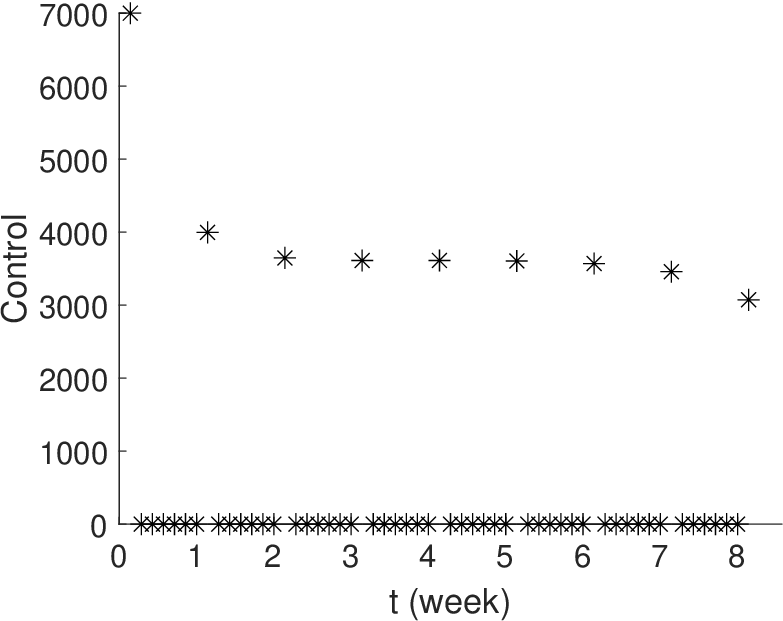} & &  \includegraphics[width=7cm]{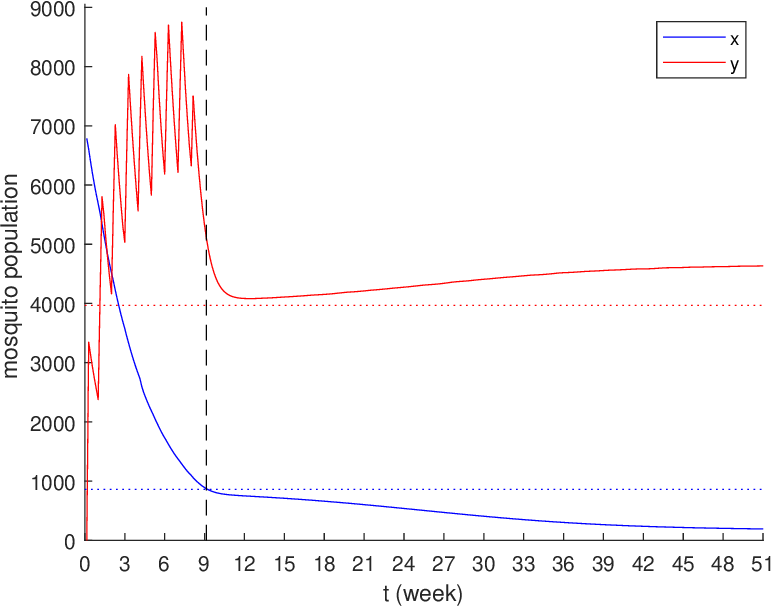} \\
    \end{tabular}
     \vspace{2mm}

    Fortnightly releases \\   \vspace{4mm}
     \begin{tabular}{ccc}
    \includegraphics[width=7cm]{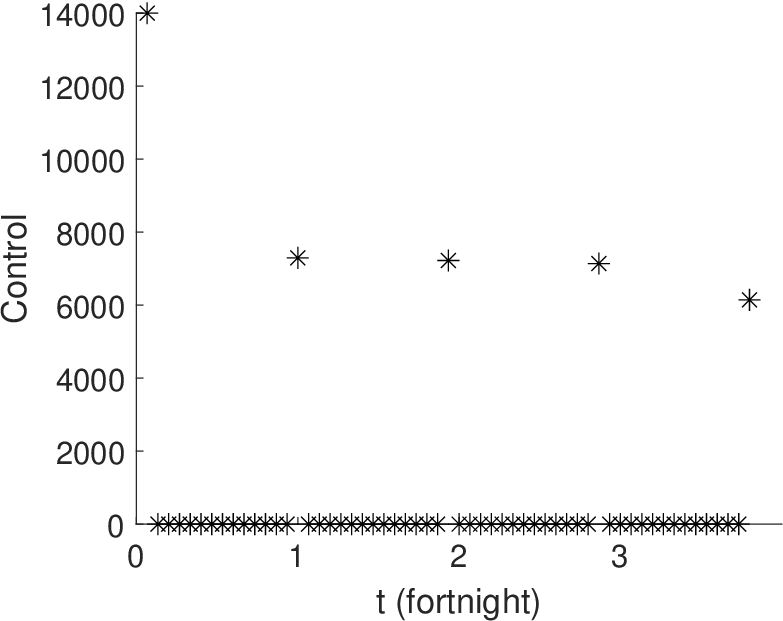} & &     \includegraphics[width=7cm]{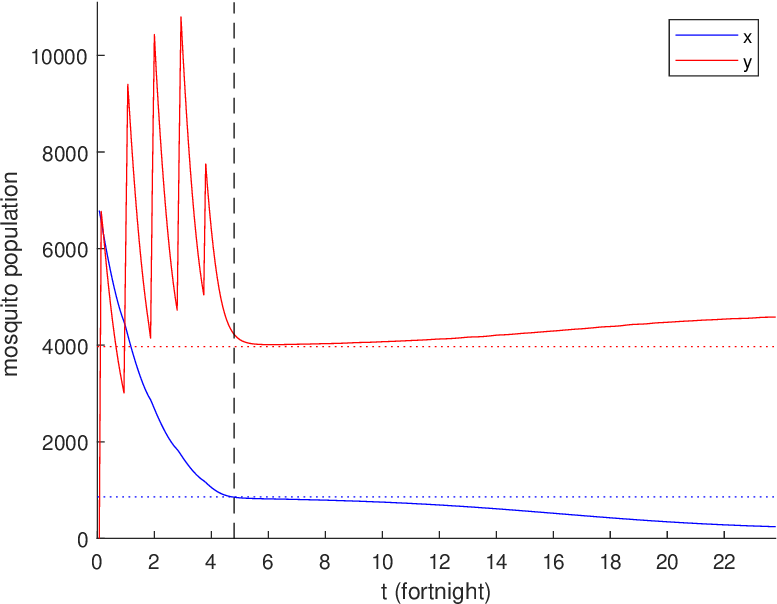} \\
    \end{tabular}
    \caption{Suboptimal impulsive control strategies for releasing \textit{w}MelPop-carrying mosquitoes \textit{(left column)} and the corresponding evolution of wild mosquitos $x(t)$ and \textit{w}MelPop-infected mosquitoes $y(t)$  \textit{(right column)}, where the vertical line marks the time of fulfillment of the condition \eqref{basin-cond}. \label{fig:wMelPop-subopt} }
\end{figure}

To simulate the suboptimal impulsive release strategies derived from the optimal continuous control $u^{*}(t)$, we construct the sequence of daily impulses $\Big\{ \widehat{U}_n^{*} \Big\}$ using the rule \eqref{subopt-imp}, and then solve the impulsive ODE system \eqref{sys-imp}. The underlying results appear in the upper rows of Figures \ref{fig:wMel-subopt} and \ref{fig:wMelPop-subopt} for the \textit{w}Mel and \textit{w}MelPop strains, respectively. The dashed vertical line in the right-column charts of Figures \ref{fig:wMel-subopt} and \ref{fig:wMelPop-subopt} marks the time of fulfillment of the condition \eqref{basin-cond}. For daily releases of \textit{w}Mel-carrying insects we have $\hat{T}^{*} =14$ days (cf. upper row of Figure \ref{fig:wMel-subopt}), while for the \textit{w}MelPop-carrying mosquitoes we have $\hat{T}^{*} =65$ days (cf. upper row of Figure \ref{fig:wMelPop-subopt}).

To propose suboptimal impulsive release strategies with more sparse frequencies, we must recall that the entire life cycle of \textit{Aedes aegypti} mosquitoes from egg to adult usually takes between one and two weeks under favorable climatic conditions. Therefore, the release frequencies must not exceed two weeks. Furthermore, when the \textit{w}Mel-infected mosquitoes are released, they will remain reproductively active for more than two weeks because their lifespan is about four weeks, and thus, they will coexist with the newborns. Therefore, we can use the rule \eqref{imp-m} to construct the weekly ($m=7$ days) and fortnightly ($m=14$ days) impulsive release sequences $\Big\{ \widehat{U}_{m,i}^{*} \Big\}$ when dealing with the \textit{w}Mel \textit{Wolbachia} strain. The sizes of suboptimal impulsive releases are marked by isolated stars in the left-hand charts of Figure \ref{fig:wMel-subopt} (middle and lower rows), and the corresponding right-hand charts display the solutions of the impulsive dynamical system \eqref{sys-imp-x}, \eqref{sys-imp-y}, \eqref{imp-mday} with $\Big\{ \widehat{U}_{7,i}^{*} \Big\}$ and $\Big\{ \widehat{U}_{14,i}^{*} \Big\},$ respectively. Note that the blue curves correspond to $x(t)$ (wild mosquito population), and the red curves stand for $y(t)$ (\textit{Wolbachia}-infected population).

On the other hand, the \textit{w}MelPop-infected mosquitoes, whose lifespan is about two weeks, will remain reproductively active for a shorter time after each release, and only a part of them may coexist with the newborns. This leads us to construct the weekly and fortnightly impulsive release sequences $\Big\{ \widetilde{U}_{m,i}^{*} \Big\}$ by the rule \eqref{imp-max} when dealing with the \textit{w}MelPop \textit{Wolbachia} strain and thus to guarantee more abundant releases. The sizes of suboptimal impulsive releases $\widetilde{U}_{m,i}^{*}$ are marked by isolated stars in the left-hand charts of Figure \ref{fig:wMelPop-subopt} for weekly (middle row, $m=7$ days) and fortnightly releases (lower row, $m=14$ days). The corresponding right-hand charts display the solutions of the impulsive dynamical system \eqref{sys-imp-x}, \eqref{sys-imp-y}, \eqref{imp-mday-max} with $\Big\{ \widetilde{U}_{7,i}^{*} \Big\}$ and $\Big\{ \widetilde{U}_{14,i}^{*} \Big\},$ respectively.

In all charts located in the right-hand columns of Figures \ref{fig:wMel-subopt} and \ref{fig:wMelPop-subopt} (and also in Figures \ref{fig:wMel_GA_graphs} and \ref{fig:wMelpop_GA_graphs} displaying the GA results), the horizontal dotted lines indicate the coordinates of the saddle-point unstable equilibrium $\mathbf{E}_u=\big( x_u, y_u \big)$. Namely, the blue dotted lines denote $x_u$, and the red dotted lines indicate $y_u$. Here it is worth pointing out that, in the case of \textit{w}Mel strain, we have that $x_u > y_u$ (see Figure \ref{fig:wMel-subopt}, right column), while for the \textit{w}MelPop strain the opposite relationship takes place, that is,  $x_u < y_u$ (see Figure \ref{fig:wMelPop-subopt}, right column). Notably, both relationships are also consonant with the position of $\mathbf{E}_u$ given in phase portraits for both \textit{Wolbachia} strains (see Figure \ref{fig:phase_planes}).

It is worthwhile to point out that all suboptimal impulsive release strategies (see the left-side charts in Figures \ref{fig:wMel-subopt} and \ref{fig:wMelPop-subopt}) are non-increasing because they naturally inherit the property of monotonicity exhibited by the continuous release strategies (see  Figures \ref{figPO1}(a) and \ref{figPO2}(a)). Thus, the monotonicity of the suboptimal impulsive release strategies allows for their easier planning from a practical standpoint.

Finally, Table \ref{tab:MPP_disc} summarizes the key indicators of the suboptimal impulsive release strategies with daily, weekly, and fortnightly frequencies. {For both \textit{Wolbachia} strains, the impulsive release strategies require releasing a more significant overall number of \textit{Wolbachia}-carrying mosquitoes than the idealistic continuous release strategies (cf. formula \eqref{int-value}), so they are in effect suboptimal.

\begin{table}[htpb]
\centering
		\caption{Key indicators of the suboptimal impulsive release strategies}
		\begin{tabular}{c|ccc|ccc} 
			\toprule
			\multirow{2}{*}{Indicators} & \multicolumn{3}{c}{\textit{w}Mel strain} & \multicolumn{3}{c}{\textit{w}MelPop strain}\\
			\cmidrule{2-4} \cmidrule{5-7}
			 & day & week & fortnight  & day & week & fortnight  \\
			\midrule
			Number of  & $\hat{T}^{*} \!=\! 14$ & $\hat{T}^{*}_7 \!=\! 2$ & $\hat{T}^{*}_{14} \!=\! 1$  & $\hat{T}^{*} \!=\! 65$ & $\hat{T}^{*}_7 \!=\! 9$ & $\hat{T}^{*}_{14} \!=\! 5$  \\
            releases &  &  &   &  &  &   \\
			Overall size & $\sum \limits_{i=1}^{14} \widehat{U}_{i}^{*} \!=\! 5966$ & $\sum \limits_{i=1}^{2} \widehat{U}_{7,i}^{*} \!=\! 5966$ & $\widehat{U}_{14,1}^{*} \!=\! 5966$ & $\sum \limits_{i=1}^{65} \widehat{U}_{i}^{*} \!=\! 33169$ & $\sum \limits_{i=1}^{9} \widetilde{U}_{7,i}^{*} \!=\! 35574$ &  $\sum \limits_{i=1}^{5} \widetilde{U}_{14,i}^{*} \!=\! 41804$ \\
            of releases &  &  &   &  &  &   \\
			\bottomrule
		\end{tabular}
		\label{tab:MPP_disc}
	\end{table}

\subsection{Release strategies yielded by the Genetic Algorithm}
\label{ga_results}

To solve the optimization problem \eqref{ga-problem}, we used MATLAB (version R2022b) to implement Algorithm \ref{alg:cap}. Here, we have set the initial value for $\epsilon$ similar to $\hat{T}^{*}$ obtained in Subsection \ref{pontr_results} for both \textit{Wolbachia} strains. Other GA-related parameters are summarized in Table \ref{tab:ga_parameters}, while the parameter values corresponding to the mosquito dynamics are taken from Table \ref{tab:model parameters}.

\begin{table}[htpb]
\centering
\caption{GA-related parameters}
\begin{tabular}{ccccc|ccc} 
\toprule
\multirow{2}{*}{Parameter} & \multirow{2}{*}{Description} &  \multicolumn{3}{c}{\textit{w}Mel strain} & \multicolumn{3}{c}{\textit{w}MelPop strain}\\
\cmidrule{3-5} \cmidrule{6-8}
  & & day & week & fortnight  & day & week & fortnight  \\
 \midrule
$N$ & population size & 100 & 100 & 100  & 100 & 100 & 100   \\
$G$ & number of generations & 100 & 100 & 100  & 100 & 100 & 100   \\
$pL$ & control upper bound & 750 & 7 $\times$ 750 & 14 $\times$ 750 & 1000 & 7  $\times$ 1000  & 14 $\times$ 1000  \\
 \bottomrule
  \end{tabular}
   \label{tab:ga_parameters}
\end{table}

It is worthwhile to point out that the GA-population size $N$ and the number of GA-generations $G$ were calibrated to obtain convergence to a globally optimal solution of the optimization problem \eqref{ga-problem}. The control upper bound $pL$, meaning the quantity of \textit{Wolbachia}-carrying insects available to release in each period of $p$ days, was chosen in accordance with \eqref{L-value}. This bound has proved again to be highly determinant for the results related to both strains since the higher (lower) its value, the lower (higher) the number releases leading to a successful control intervention.

As indicated in Step \textbf{1} of Algorithm \ref{alg:cap}, $T$ is a multiple of $p=1,7,14$, and in the sequel, we denote by $T_p^{*}$ the optimal time rendered by Genetic Algorithm. The objective function \eqref{functional} directly expresses the overall number of \textit{Wolbachia}-carrying mosquitoes to be released during the intervention. Figures \ref{fig:wMel_GA_graphs} and \ref{fig:wMelpop_GA_graphs} display the graphs of the best GA-based control strategies \textit{(left columns)} and the associated system behavior \textit{(right columns)} for the \textit{w}Mel and \textit{w}MelPop \textit{Wolbachia} strains, respectively.

In both figures, the daily release strategies ($p=1$) are presented in the upper rows, weekly release strategies ($p=7$) are given in the middle rows, and fortnightly release strategies ($p=14$) appear in the lower rows. The optimal time $T_p^{*}, p=1,7,14$ and the total number of \textit{Wolbachia}-carrying mosquitoes needed for the successful intervention, $J \big( u^{*}, T^{*}_p \big)$ are given in Table \ref{tab:ga_results}. In this table and also in Figures \ref{fig:wMel_GA_graphs} and \ref{fig:wMelpop_GA_graphs} we observe again that $T^{*}_p, p=1,7,14$ is far shorter for the \textit{w}Mel strain (about $1-2$ weeks) compared to the \textit{w}MelPop strain (around $9$ weeks).

\begin{table}[htpb]\centering
\caption{Key indicators of the release strategies yielded by Genetic Algorithm}
\begin{tabular}{c|ccc|ccc} 
\toprule
\multirow{2}{*}{Indicators}  & \multicolumn{3}{c}{\textit{w}Mel strain} & \multicolumn{3}{c}{\textit{w}MelPop strain}\\
\cmidrule{2-4} \cmidrule{5-7}
   & day & week & fortnight  & day & week & fortnight  \\
 \midrule
Number of releases, $T^{*}_p$ & 11 & 2 & 1  & 60 & 9 & 5  \\
Overall size of releases $J \big( u^{*},T^{*}_p \big)$ & 5436 & 5226 & 4956 & 24481 & 27259 &  31323 \\
 \bottomrule
  \end{tabular}
   \label{tab:ga_results}
\end{table}

As shown on the left upper charts of Figures \ref{fig:wMel_GA_graphs} and \ref{fig:wMelpop_GA_graphs} as well as reflected in Table \ref{tab:ga_results}, the total number of effective daily releases $(p=1)$ rendered by GA for both strains is a bit smaller than in the case of suboptimal impulsive strategies for daily releases presented in Subsection \ref{pontr_results} (cf. left-hand upper charts in Figures \ref{fig:wMel-subopt} and \ref{fig:wMelPop-subopt} and the data from Table \ref{tab:MPP_disc}). Indeed, in the case of \textit{w}Mel strain intervention, the trajectories $x(t)$ and $y(t)$ of the dynamical system \eqref{system} with GA-produced $u(t)$ will fulfill the conditions \eqref{basin-cond} after $11$ daily releases. In contrast, $14$ daily impulsive releases will be needed to move the trajectories $x(t)$ and $y(t)$ of \eqref{sys-imp} in the attraction basin of $\mathbf{E}_s$,  the equilibrium of stable coexistence. The case of \textit{w}MelPop strain intervention exhibits a similar outcome: $60$ daily releases when GA-based release strategies are employed versus $65$ daily releases using suboptimal impulsive strategies.

Let us now address the total number of \textit{Wolbachia}-carrying mosquitoes needed for successful intervention, another critical indicator of the release strategies. As shown in Table \ref{tab:ga_results}, this number is reduced as $p$ grows from day to a week and further to fortnight in the case of the \textit{w}Mel strain, and the situation is opposite in the case of the \textit{w}MelPop strain. Furthermore, by comparing the data from Tables  \ref{tab:MPP_disc} and \ref{tab:ga_results}, it can be observed that the GA-based release strategies outperform the suboptimal impulsive release strategies, derived from the Pontryagin Maximum Principle, for both \textit{Wolbachia} strains and for all release frequencies. Indeed, the GA-based release strategies, compared to the impulsive ones, reduce the overall quantity of \textit{Wolbachia}-carrying insects, depending on the release frequencies,  by 9-17\% (in the case of \textit{w}Mel strain) and 23-26\% (in the case of \textit{w}MelPop strain).

\begin{figure}[H]
    \centering
    Daily releases: $p=1$ \\ \vspace{4mm}
        \begin{tabular}{ccc}
    \includegraphics[width=7cm]{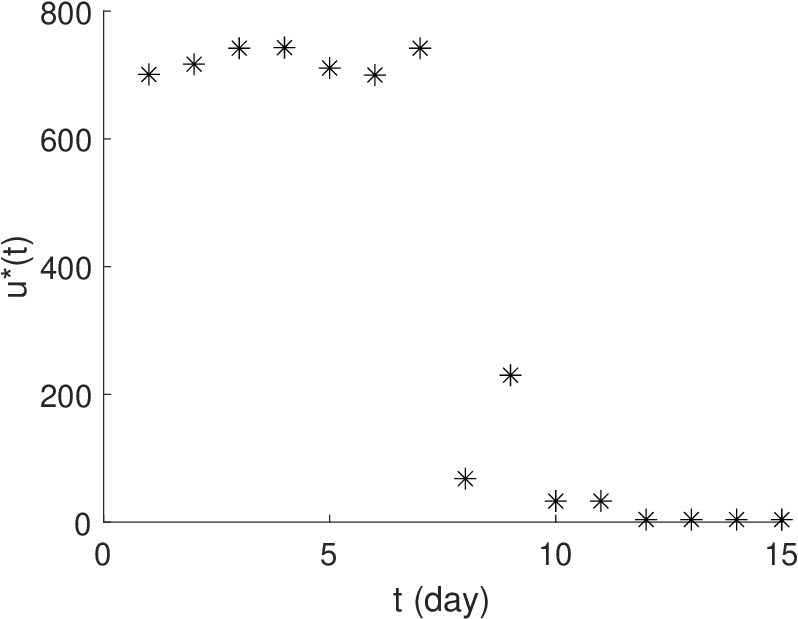} & & \includegraphics[width=7cm]{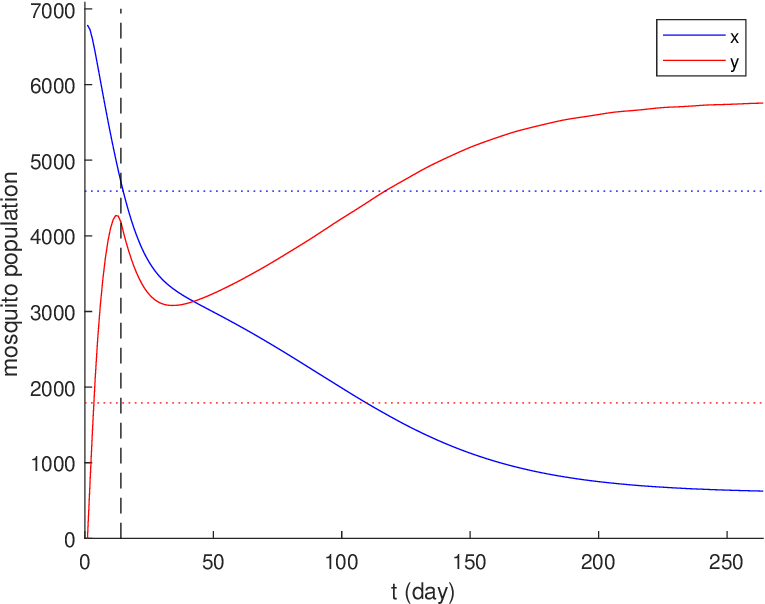} \\
    \end{tabular}
    \vspace{2mm}

    Weekly releases: $p=7$ \\ \vspace{4mm}
         \begin{tabular}{ccc}
   \includegraphics[width=7cm]{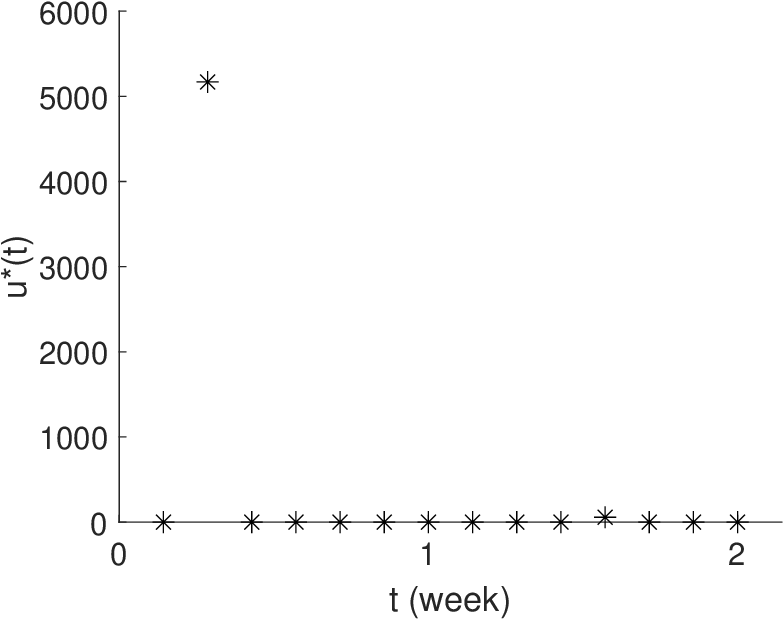} & &  \includegraphics[width=7cm]{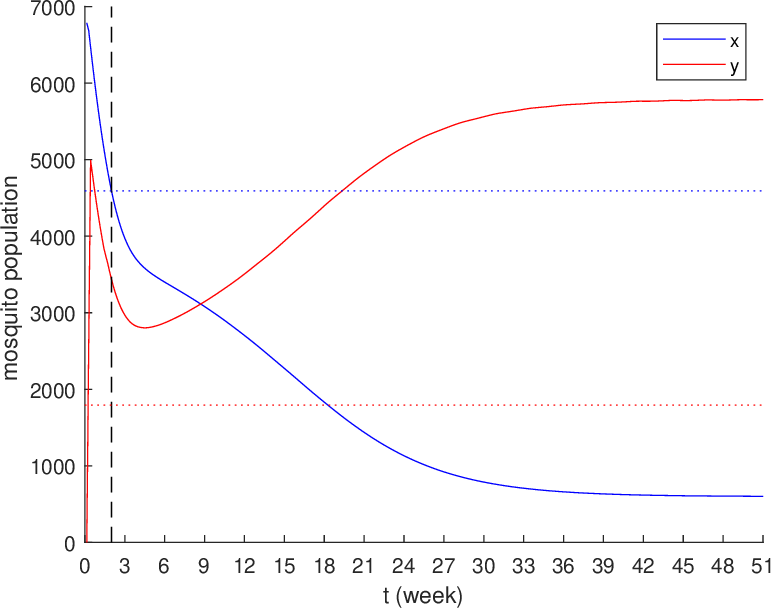}  \\
    \end{tabular}
     \vspace{2mm}

    Fortnightly releases: $p=14$ \\   \vspace{4mm}
     \begin{tabular}{ccc}
    \includegraphics[width=7cm]{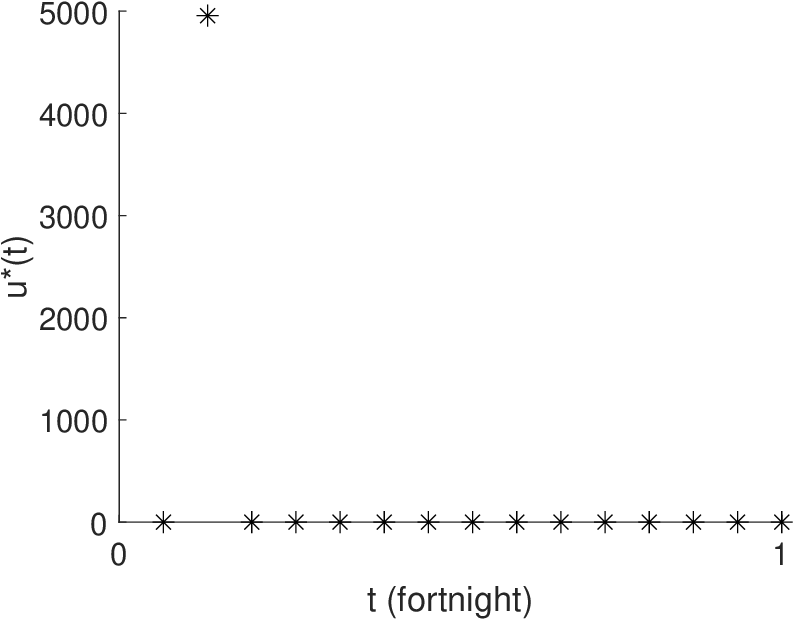} & &  \includegraphics[width=7cm]{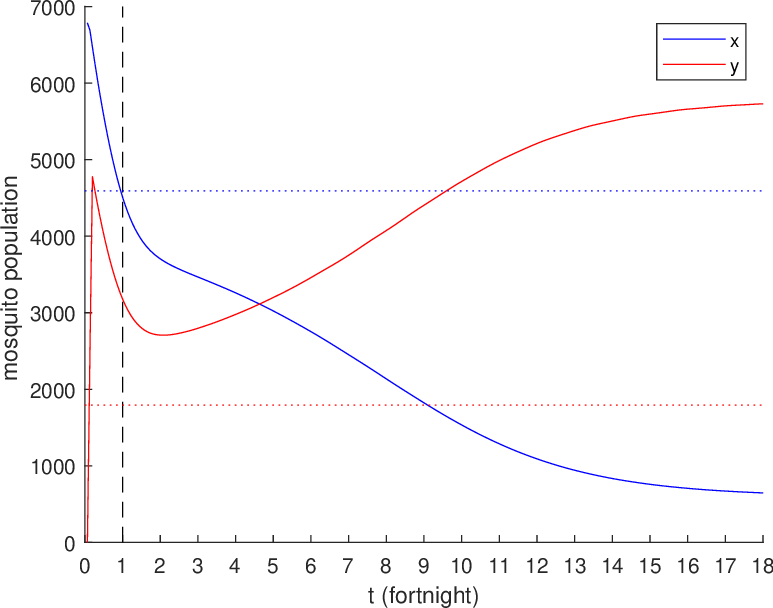} \\
    \end{tabular}
    \caption{GA-optimal strategies for releasing \textit{w}Mel-carrying mosquitoes \textit{(left column)} and the corresponding evolution of wild mosquitos $x(t)$ and \textit{w}Mel-infected mosquitoes $y(t)$  \textit{(right column)}, where the vertical line marks the time of fulfillment of the conditions \eqref{basin-cond}. \label{fig:wMel_GA_graphs}}
\end{figure}

\begin{figure}[H]
    \centering
    Daily releases: $p=1$ \\ \vspace{4mm}
        \begin{tabular}{ccc}
    \includegraphics[width=7cm]{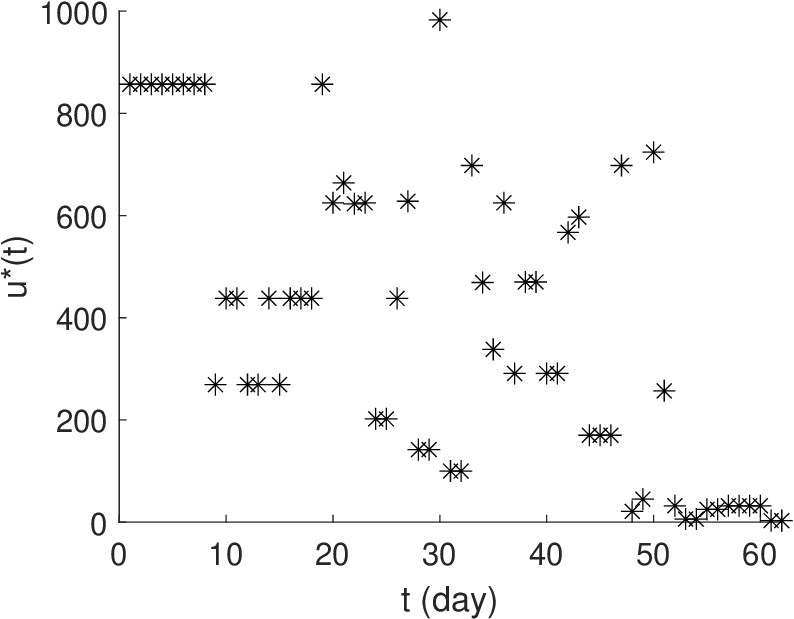} & & \includegraphics[width=7cm]{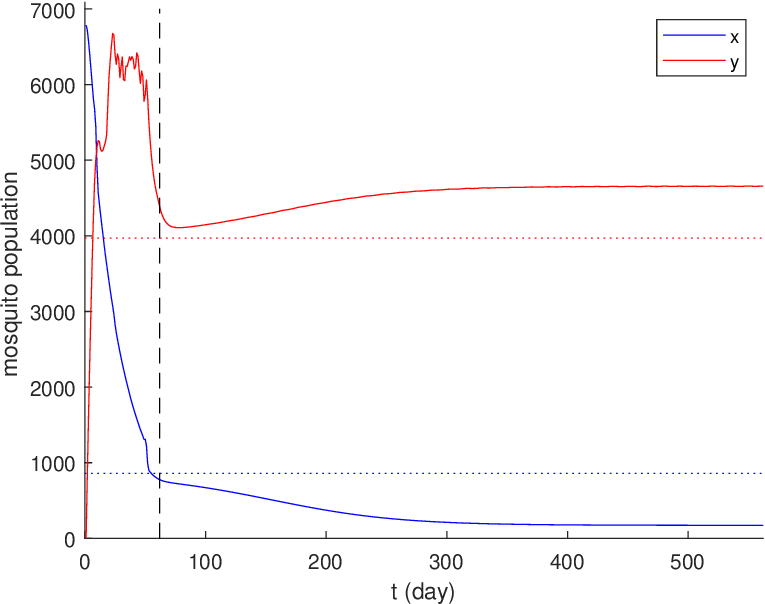} \\
    \end{tabular}
    \vspace{2mm}

    Weekly releases: $p=7$ \\ \vspace{4mm}
         \begin{tabular}{ccc}
   \includegraphics[width=7cm]{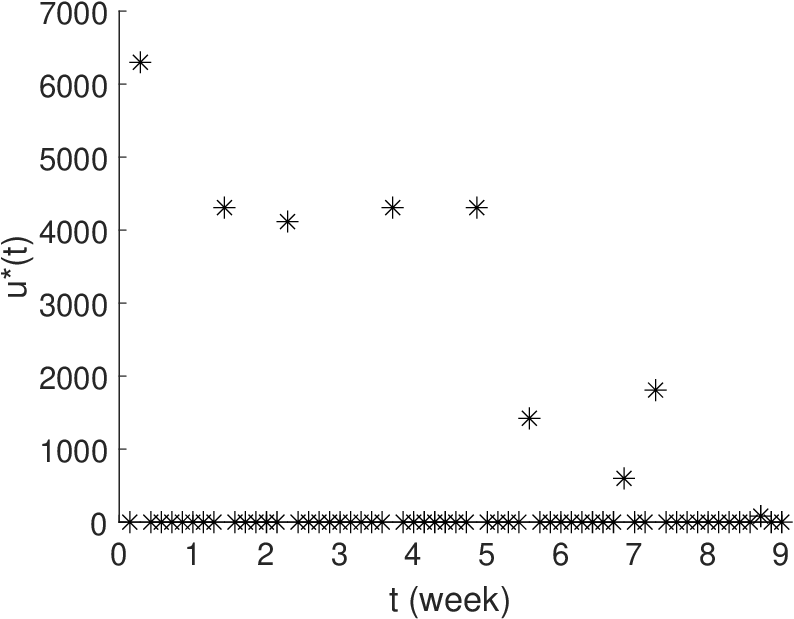} & &  \includegraphics[width=7cm]{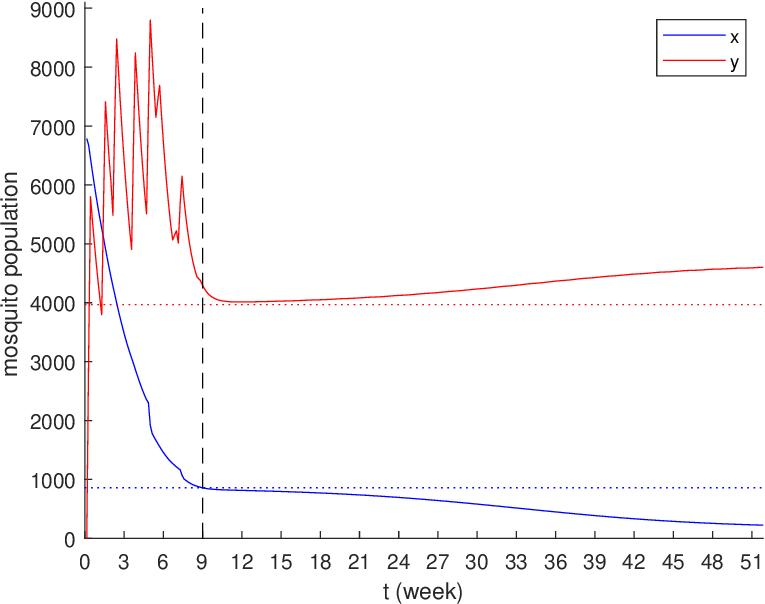}  \\
    \end{tabular}
     \vspace{2mm}

    Fortnightly releases: $p=14$ \\   \vspace{4mm}
     \begin{tabular}{ccc}
    \includegraphics[width=7cm]{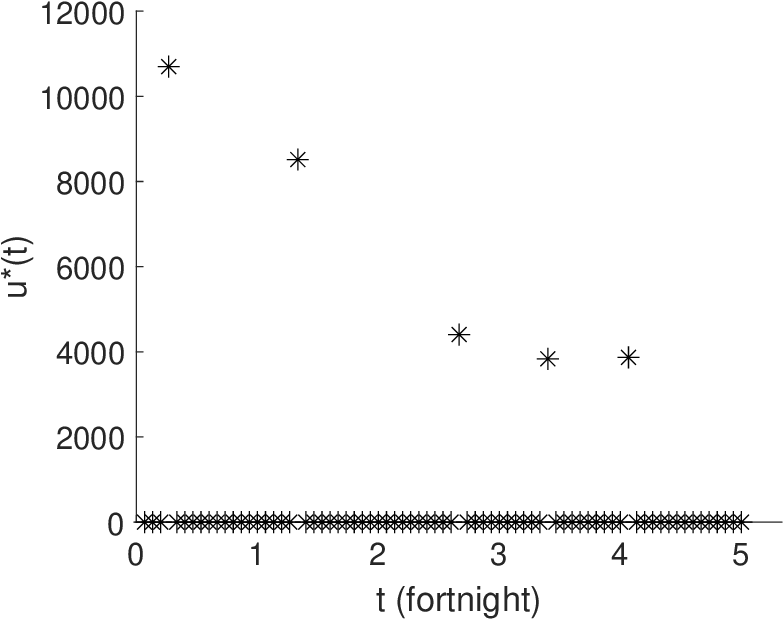} & &  \includegraphics[width=7cm]{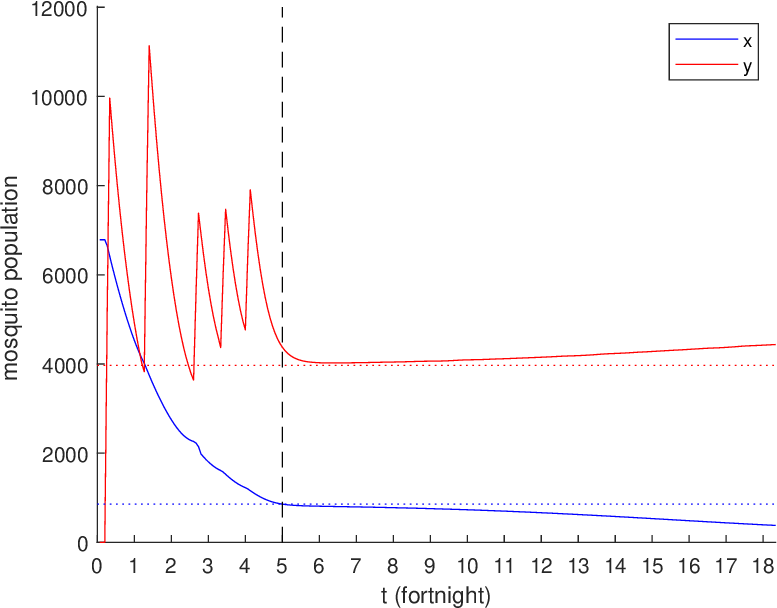} \\
    \end{tabular}
    \caption{GA-optimal strategies for releasing \textit{w}MelPop-carrying mosquitoes \textit{(left column)} and the corresponding evolution of wild mosquitos $x(t)$ and \textit{w}Mel-infected mosquitoes $y(t)$  \textit{(right column)}, where the vertical line marks the  time of fulfillment of the conditions \eqref{basin-cond}. \label{fig:wMelpop_GA_graphs}}
\end{figure}

Another essential feature of the GA-generated release strategies is their non-monotonicity, especially when $p=1$ (daily releases, cf. upper left-hand charts in Figures \ref{fig:wMel_GA_graphs} and \ref{fig:wMelpop_GA_graphs}). This marks their structural difference from the suboptimal impulsive release strategies presented in Subsection \ref{pontr_results} and implies more careful planning of the releases in the practical settings because the releases require previous mass-rearing of the \textit{Wolbachia}-carrying insects. However, more sparse  release strategies rendered by the genetic algorithm ($p=7$ and $p=14$, see the left-hand charts in the middle and lower rows of Figures \ref{fig:wMel_GA_graphs} and \ref{fig:wMelpop_GA_graphs} are more feasible for the ahead planning practical implementation. Furthermore, it is worthwhile to point out that the GA routine given by Algorithm \ref{alg:cap} also optimizes the release day within $p$ days when dealing with sparser releases ($p=7,14$); see the left-hand charts in the middle and lower rows of Figures \ref{fig:wMel_GA_graphs} and \ref{fig:wMelpop_GA_graphs}. In contrast, the suboptimal impulsive release strategies anticipate only ``equidistant'' releases, e.g., every Monday ($m=7$) or every other Monday ($m=14$), as displayed by the left-hand charts in the middle and lower rows of Figures \ref{fig:wMel-subopt} and \ref{fig:wMelPop-subopt}.

\section{Discussion and concluding remarks}
\label{discussion}

The numerical results presented in Section \ref{results} illustrated two rational and workable approaches for implementing \textit{Wolbachia}-based biocontrol in practical settings. Both approaches facilitate establishing one of the two \textit{Wolbachia} strains (\textit{w}Mel or \textit{w}MelPop) in the target geographic area by periodic releases of \textit{Wolbachia}-infected mosquitoes. However, the release strategies produced by both approaches have resulted in structurally different, and their core indicators also present some dissimilarities.

To understand why these optimization approaches have produced dissimilar release strategies, we should revisit the mathematical backgrounds of the optimization problems laid into the basis of each optimization approach. Here we recall first that the sequences of suboptimal impulsive releases $\big\{ \widehat{U}^{*}_{n} \big\}, \big\{ \widehat{U}^{*}_{m,i} \big\}$ and  $\big\{ \widetilde{U}^{*}_{m,i} \big\}$ (see formulas \eqref{subopt-imp}, \eqref{imp-m}, \eqref{imp-max}) were constructed from the optimal control $u^{*}(t), t \in \big[ 0, T^{*} \big]$, which is a continuous real function. This function was an optimal solution to the OCP \eqref{ocp} where the class of admissible control functions consisted of real piece-wise continuous functions satisfying the condition \eqref{umax}. Thus, the optimization was first performed over the class of real piece-wise continuous functions, and the optimal solution then gave rise to a sequence of suboptimal impulses.

On the other hand, the sequence of near-optimal decisions $\big( u^*(0), u^*(1), \ldots, u^*(T^{*}) \big)$ produced by the genetic algorithm was the result of direct optimization performed over the class of discrete sequences, which is a different set of admissible controls. Furthermore, the GA routine given by Algorithm \ref{alg:cap} also allowed the optimization of the release day within a week or fortnight period, while the suboptimal impulsive sequences only anticipated equidistant releases. This explains why the GA-produced release strategies are structurally dissimilar to the impulsive release strategies derived from the maximum principle and display better quantitative outcomes.

In contrast to recent works \cite{Bliman2023,Vicencio2023} where constant-size periodic releases were suggested for \textit{Wolbachia}-based biocontrol, the impulsive release strategies derived presented in Subsection \ref{pontr_results} bear variable release sizes that are monotone decreasing and thus require a smaller total number of \textit{Wolbachia}-carrying insects for successful intervention. Furthermore, the GA-yielded release strategies derived in Subsection \ref{ga_results} virtually anticipate non-monotone variable release sizes that apparently help reduce even more the overall quantity of \textit{Wolbachia}-carriers needed for intervention.

Contrasting the data from Tables \ref{tab:MPP_disc} and \ref{tab:ga_results}, we observe that even though both optimization approaches rendered almost the same number of releases, the outcomes of the GA-yielded release strategies are slightly better for both \textit{Wolbachia} strains in the case of daily releases. Moreover, the overall number of \textit{Wolbachia}-carrying insects required for successful intervention is always less in the case of the GA approach. From a practical standpoint, establishing \textit{Wolbachia} in the wild population of \textit{Aedes aegypti} mosquitoes by releasing a smaller number of \textit{Wolbachia}-carrying insects is beneficial for the following reasons:
\begin{itemize}
  \item
  Smaller release sizes will produce lower costs related to mass-rearing of \textit{Wolbachia}-infected mosquitoes, such as infrastructure, lab equipment, personnel, nutrients, etc.
  \item
  Releasing fewer females will entail a lower nuisance among the human inhabitants because  \textit{Wolbachia}-infected female mosquitoes still need to bite people to mature their eggs.
  \item
  Even though \textit{Wolbachia} induces the virus blockage in mosquitoes, a share of \textit{Wolbachia}-infected female mosquitoes is still capable of transmitting the viruses during biting; therefore, releasing a smaller number of females may potentially reduce the number of human viral infections.
\end{itemize}

Regarding the frequency of releases, one should keep in mind that, besides mass-rearing, each release also entails significant costs related to logistics and personnel for sorting, transporting, and releasing the insects. In this sense, our simulations indicate that fortnightly release strategies yielded by GA are a cheaper option for establishing \textit{Wolbachia} in some local areas because they require releasing not only a smaller overall quantity of mosquitoes but also fewer releases. However, taking into account that about $20\%$ of \textit{w}Mel-infected female mosquitoes and around $5\%$ of \textit{w}MelPop-infected ones are still capable of transmitting the virus to humans \cite{Dorigatti2018}, the more abundant fortnightly releases may drastically increase the vectorial density in the areas of releases, cause more nuisance to local residents, and temporarily enhance the number of human infections. Thus, our simulation results may help the local healthcare authorities to revise all such arguments before planning an intervention based on either \textit{w}Mel or \textit{w}MelPop strain.

Finally, based on the results presented in Section \ref{results}, we recommend GA-based optimization for planning the release programs since this approach performs optimization directly over the class of discrete sequences of decisions. Notably, the GA approach can be adapted to a desired frequency of the releases and different initial sizes of the wild mosquito population (note that all simulations in Section \ref{results} correspond to the wild population inhabiting 1 hectare). Furthermore, this approach can also be tested for more sophisticated population dynamics of wild and \textit{Wolbachia}-carrying mosquitoes, including stage structure, sex structure, or spatial distribution. For such tests, the dynamical system \eqref{system} should be replaced by a more complex one, and the constraints \eqref{x_threshold}, \eqref{y_threshold} must be adjusted accordingly.

\section*{Acknowledgments}

The authors thank the Brazilian foundation FAPESP (Grant Number 2020/10964-0), CNPq (Grant Number 306518/2022-8), CAPES (Finance Code 001), and the Colombian Ministry of Science, Technology, and Innovation --- \textit{MINCIENCIAS} (Contract No. 80740-224-2021) for financial support.



\bibliographystyle{elsarticle-num}
\bibliography{WB-opt-bib}







\end{document}